\documentclass[3pt]{article}

\usepackage{graphicx}\usepackage{amssymb}
\usepackage{amsmath}\usepackage{amsfonts}
\numberwithin{equation}{section}

\begin{document}

\title{The Multilevel Finite Element Discretizations Based on Local Defect-Correction
for Nonsymmetric Eigenvalue Problems}
\author{ { Yidu Yang, Jiayu Han} \\
{\small School of Mathematics and Computer Science, }\\{\small
Guizhou Normal University,  Guiyang,  $550001$,  China} }\date{~}
\pagestyle{plain} \textwidth 145mm \textheight 215mm \topmargin 0pt
\maketitle

\indent{\bf\small Abstract~:~} {\small Based on the work of Xu and Zhou [Math.Comput., 69(2000), pp.
881-909], we establish  new three-level and multilevel finite
element discretizations by local defect-correction technique. Theoretical analysis and numerical experiments show that the schemes
are simple and easy to carry out, and can be used to solve singular
nonsymmetric eigenvalue problems efficiently. We also discuss the
local error estimates of finite element approximations; it's a new
feature here that the estimates apply to the local domains
containing corner points.
}\\
\indent{\bf\small Keywords~:~\scriptsize} {\small nonsymmetric eigenvalue problems, finite element, multilevel
discretization, local refinement, local defect-correction.
\\\indent{\bf\small 1991 MSC~:~\scriptsize}    code 65N25,65N30}

\section{Introduction}

\indent Nonsymmetric elliptic eigenvalue problems have important
physical background, such as convection-diffusion in fluid
mechanics, environmental problems and so on. Thus, finite element
methods for solving this problem become an important topic which has
attracted the attention of mathematical and physical fields:
\cite{bramble} discussed  a priori error estimates,
\cite{carstensen,han,gedicke,heuveline1,heuveline2,yang1}  a
posteriori error estimates and adaptive algorithms, \cite{naga}
function value recovery algorithms, \cite{kolman,yang2} two level
algorithms, \cite{lu,yang3}  extrapolation methods,
\cite{carstensen} an adaptive homotopy approach, etc. This paper turns to
discuss finite element multilevel discretization based on
local defect-correction.\\
\indent For elliptic boundary value problem, Xu and Zhou \cite{xu1}
combined two-grid finite element discretization scheme with the
local defect-correction technique to propose a general and powerful
parallel-computing technique. This technique has been used and
developed by many scholars, for instance, successfully applied to
Stokes equation~(see \cite{he1,he2}), especially, Xu and Zhou
\cite{xu3}, Dai and Zhou \cite{dai}, Bi and Yang etc \cite{bi}
developed this method and established local and parallel three-level
finite element discretizations for symmetric elliptic singular
eigenvalue problems
(including the electronic structure problems).\\
\indent In this paper, we further apply local defect-correction
technique proposed by Xu and Zhou  to~nonsymmetric elliptic singular
eigenvalue problems, our work has the following features. (1) We
extend local and parallel three-level finite element discretizations
for symmetric eigenvalue problems established by Dai and Zhou
\cite{dai} to~nonsymmetric eigenvalue problems. (2)~Based on
\cite{bi}, we establishes new multilevel finite element
discretization by local refinement, this scheme repeatedly makes
defect correction on finer and finer local meshes to make up for
abrupt changes of local mesh size caused by three level scheme. And
theoretical analysis and numerical experiments show that our schemes
are simple and easy to carry out, and can be used to solve singular
nonsymmetric eigenvalue problems. Numerical experiments show that,
compared with the adaptive homotopy approach in \cite{carstensen},
our algorithm seems to be more efficient. (3) For the nonsymmetric
problems, based on the work of \cite{wahlbin,xu1}, we discuss the
local error estimates of  finite element approximations; it’s a new
feature here that the estimates apply to the local domains
containing corner points, see Lemmas 2.3-2.4 and
Remark 2.2 in this paper.\\

\indent In this paper, regarding the basic theory of  finite
elements, we refer to
\cite{babuska,brenner,ciarlet,oden}.\\

\section {preliminaries}

\indent Consider the nonsymmetric elliptic differential operator
eigenvalue problem:
\begin{eqnarray}\label{s2.1}
&&Lu\equiv -\sum\limits_{i,j=1}^{d}\partial_{j}(a_{ij}(x) \partial_{
i}u) +\sum\limits_{i=1}^{d}b_{i}(x)\partial_{i}u +c(x)u=\lambda m(x)
u,~~~\mathrm{in}~\Omega,~~~~~~\\\label{s2.2} &&~~~~~~~~~u=0,~~~\mathrm{on}~
\partial\Omega,
\end{eqnarray}
where  $\Omega\subset R^{d}$, $d\geq 2$, be a polyhedral bounded
domain with boundary $\partial\Omega$, $\partial_{
i}u=\frac{\partial u}{\partial x_{i}}$, $i=1,2,\cdots,d$.\\

\indent Let
\begin{eqnarray*}
&&a(u,v)=\int\limits_{\Omega}(\sum\limits_{i,j=1}^{d}a_{ij}\partial_{i}u\partial_{j}\overline{v}
+\sum\limits_{i=1}^{d}b_{i}\partial_{i}u\overline{v}
+cu\overline{v})\mathrm{d}x,\\
&&b(u,v)=\int\limits_{\Omega} m u\overline{v}\mathrm{d}x.
\end{eqnarray*}

\indent The variational form associated with
(\ref{s2.1})-(\ref{s2.2}) is given by: find $\lambda \in
\mathcal{C}$, $u\in H_{0}^{1}(\Omega),
~\|u\|_{0}=1$, satisfying
\begin{eqnarray}\label{s2.3}
a(u,v)=\lambda b(u,v),~~~\forall~ v\in H_{0}^{1}(\Omega).
\end{eqnarray}

Assume that $a_{i,j}, b_{i}\in W_{1, \infty}(\Omega)$, $c\in L_{
\infty}(\Omega)$ are given real or complex functions on $\Omega$,
$m\in L_{ \infty}(\Omega)$ is a given real function which is bounded
below by a positive constant on $\Omega$. $L$ is assumed to be
uniformly strongly elliptic in $\Omega$, i.e., there is a positive
constant $a_{0}$ such that
\begin{eqnarray}\label{s2.4}
Re \sum\limits_{i,j=1}^{d}a_{i,j}\xi_{i}\xi_{j}\geq
a_{0}\sum\limits_{i=1}^{d}\xi_{i}^{2},~~~\forall x\in\Omega,~\forall
(\xi_{1},\xi_{2},\cdots,\xi_{d})\in R^{d}.
\end{eqnarray}
Let $b=\max\limits_{1\leq i\leq d; x\in \Omega}|b_{i}(x)|$. we
assume without loss that $Re~c \geq a_{0}/2 +b^{2}/(2a_{0})$ since
adding a constant$\times m(x)$ to $c(x)$ only shifts the
eigenvalues. Under above assumptions, we have
\begin{eqnarray}\label{s2.5}
Re ~a(u,u)\geq \frac{1}{2}a_{0}\|u\|_{1}^{2},~~~\forall u\in
H^{1}(\Omega);
\end{eqnarray}
and there are constants $M_{1}$ and $M_{2}$ such that
\begin{eqnarray}\label{s2.6}
&&|a(u,v)|\leq M_{1}\|u\|_{1}\|v\|_{1},~~~\forall u, v\in
H^{1}(\Omega),\\\label{s2.7} &&|b(u,v)|\leq
M_{2}\|u\|_{0}\|v\|_{0},~~~\forall u, v\in L_{2}(\Omega).
\end{eqnarray}

For $D\subset\Omega_{0}\subset \Omega$, we use
$D\subset\subset\Omega_{0}$ to mean that $dist(\partial D\backslash\partial\Omega, \partial\Omega_{0}\backslash\partial\Omega)>0$.\\
\indent Assume that $\pi_{h}(\Omega)=\{\tau\}$ is a mesh of $\Omega$
with mesh-size function $h(x)$ whose value is the diameter
$h_{\tau}$ of the element $\tau$ containing $x$, and
$h(\Omega)=\max\limits_{x\in\Omega}h(x)$ is the mesh diameter of
$\pi_{h}(\Omega)$. We write $h(\Omega)$ as $h$ for simplicity. Let
$V_{h}(\Omega)\subset C(\overline{\Omega})$, defined on
$\pi_{h}(\Omega)$, be a space of piecewise polynomials, and
$V_{h}^{0}(\Omega)=V_{h}(\Omega)\cap H_{0}^{1}(\Omega)$. Given
$G\subset\Omega$, we define $\pi_{h}(G)$ and $V_{h}(G)$  to be the
restriction of $\pi_{h}(\Omega)$ and $V_{h}(\Omega)$ to $G$,
respectively, and
$$V_{h}^{0}(G)=V_{h}(G)\cap H_{0}^{1}(G),~~~
V_{0}^{h}(G)=\{v\in V_{h}^{0}(\Omega): \mathrm{supp}~v\subset\subset G\}.$$
For any $G\subset \Omega$ mentioned in this paper, we assume that it
aligns with $\pi_{h}(\Omega)$ when necessary.\\
\indent In this paper, $C$ denotes a positive constant independent
of $h$, which may not be the same constant in different places. For
simplicity, we use the symbol $a\lesssim b$ to mean that $a\leq C
b$.\\

\indent We adopt the following assumptions in \cite{xu1} for meshes
and finite
element space.\\
\indent (A0)~ There exists $\nu\geq 1$ such that $
h(\Omega)^{\nu}\lesssim h(x),~~ \forall x\in \Omega$.\\
 \indent (A1)~
There exists $r\geq 1$ such that for $w\in H_{0}^{1}(\Omega)\cap
H^{1+t}(\Omega)$,
\begin{eqnarray*}
\inf\limits_{v\in
V_{h}^{0}(\Omega)}(\|h^{-1}(w-v)\|_{0}+\|w-v\|_{1})\lesssim
h^{t}\|w\|_{1+t},~~~0\leq t\leq r.
\end{eqnarray*}
\indent (A2)~ {\em Inverse Estimate.} For any $v\in
V_{h}(\Omega_{0})$, $ \|v\|_{1,\Omega_{0}}\lesssim
\|h^{-1}v\|_{0,\Omega_{0}}$.\\
 \indent (A3)~{\em
Superapproximation.} For $G\subset \Omega_{0}$, let $\omega\in
C^{\infty}(\bar{\Omega})$ with $\mathrm{supp}~\omega \subset\subset
G$,
then for any $w\in V_{h}(G)$, $w|_{\partial G\cap\partial\Omega}=0$,
there exists $v\in V_{0}^{h}(G)$ such that
$$ \|h^{-1}(\omega w-v)\|_{1,G}\lesssim \|w\|_{1,G}.$$

\indent Let $\pi_{h}(\Omega)$ consist of shape-regular simplices and
(A0) hold, and let $V_{h}(\Omega)\subset C(\overline{\Omega})$ be a
space of piecewise polynomials of degree $\leq r$ defined on $\pi_{h}(\Omega)$, then from \cite{xu1} we know that (A1)-(A3) are valid for this $V_{h}(\Omega)$.\\

\indent The finite element approximation of (\ref{s2.3}) is given
by: find $\lambda_{h}\in C$, $u_{h} \in V_{h}^{0}(\Omega),~
\|u_{h}\|_{0}=1$, satisfying
\begin{eqnarray}\label{s2.8}
a(u_{h},v)=\lambda_{h} b(u_{h},v),~~~ \forall v \in
V_{h}^{0}(\Omega).
\end{eqnarray}

\indent Thanks to \cite{babuska}, we know the adjoint problem of
(\ref{s2.1})-(\ref{s2.2}) is:
\begin{eqnarray}\label{s2.9}
&&L^{*}u^{*}\equiv
 -\sum\limits_{i,j=1}^{d}\partial_{i}(\overline{a}_{ij}\partial_{j} u^{*})-
\sum\limits_{i=1}^{d}\partial_{i}(\overline{b}_{i}u^{*})
+\overline{c}u^{*}=\lambda^{*}m u^{*},~~~in\Omega,\\\label{s2.10}
&&~~~~~~~~~~~u^{*}=0,~~~on \partial\Omega.
\end{eqnarray}
The corresponding variational form and discrete variational form of
(\ref{s2.9})-(\ref{s2.10}) are given by: find $\lambda^{*} \in
\mathcal{C}$, $u^{*}\in H_{0}^{1}(\Omega), \|u^{*}\|_{0}=1$,
satisfying
\begin{eqnarray}\label{s2.11}
a(v,u^{*})=\overline{\lambda^{*}} b(v,u^{*}),~~~~\forall~ v\in
H_{0}^{1}(\Omega);
\end{eqnarray}
find $\lambda_{h}^{*} \in \mathcal{C}$,$u_{h}^{*}\in V_{h},
\|u_{h}^{*}\|_{0}=1$, satisfying
\begin{eqnarray}\label{s2.12}
a(v,u_{h}^{*})=\overline{\lambda_{h}^{*}}
b(v,u_{h}^{*}),~~~~\forall~ v\in V_{h}^{0}(\Omega).
\end{eqnarray}
Note that the primal and dual eigenvalues are connected via
$\lambda=\overline{\lambda^{*}}$ and
$\lambda_{h}=\overline{\lambda_{h}^{*}}$.\\

\indent Throughout this paper, we will assume that
(\ref{s2.5})-(\ref{s2.7}) hold. Thus from Lax-Milgram theorem
we know the source problem associated with (\ref{s2.3}) and (\ref{s2.11}) admits  an unique solution, respectively.
The discrete source problem associated (\ref{s2.8}) and
(\ref{s2.12}) admits an unique solution, respectively.\\

 \indent Define the
solution operator $T:L_{2}(\Omega)\to H_{0}^{1}(\Omega)$ and $T_{h}:
L_{2}(\Omega)\to V_{h}^{0}(\Omega)$ as follows:
\begin{eqnarray}\label{s2.13}
a(Tg,v)&=&b(g,v),~~~ \forall v \in H_{0}^{1}(\Omega),\\\label{s2.14}
a(T_{h}g,v)&=&b(g,v),~~~ \forall v \in V_{h}^{0}(\Omega).
\end{eqnarray}

And (\ref{s2.3}) and (\ref{s2.8}) have the equivalent operator form
(\ref{s2.15}) and (\ref{s2.16}), respectively.
\begin{eqnarray}\label{s2.15}
&&Tu=\lambda^{-1} u,\\\label{s2.16}
&&T_{h}u_{h}=\lambda_{h}^{-1}u_{h}.
\end{eqnarray}

Define the solution operator $T^{*}: L_{2}(\Omega)\to
H_{0}^{1}(\Omega)$ and $T_{h}^{*}: L_{2}(\Omega)\to
V_{h}^{0}(\Omega)$ satisfying
\begin{eqnarray}\label{s2.17}
&&a(v, T^{*}f)=b(v,f),~~~\forall~v\in
H_{0}^{1}(\Omega),\\\label{s2.18} &&a(v,
T_{h}^{*}f)=b(v,f),~~~\forall~v\in V_{h}^{0}(\Omega).
\end{eqnarray}
And (\ref{s2.11}) and (\ref{s2.12}) have the equivalent operator
forms (\ref{s2.19}) and (\ref{s2.20}), respectively.
\begin{eqnarray}\label{s2.19}
&&T^{*}u^{*}=\lambda^{*-1} u^{*},\\\label{s2.20}
&&T_{h}^{*}u_{h}^{*}=\lambda_{h}^{*-1}u_{h}^{*}.
\end{eqnarray}

\indent It can be proved that $T$ is completely continuous, and
$T^{*}$ is the adjoint operator of $T$ in the sense of inner product
$b(\cdot,\cdot)$. In fact,
\begin{eqnarray*}
&&b(Tu, v)=a(Tu, T^{*}v)=b(u, T^{*}v),~~~\forall u,v\in L_{2}(\Omega),\\
&&b(T_{h}u, v)=a(T_{h}u, T_{h}^{*}v)=b(u, T_{h}^{*}v),~~~\forall
u,v\in L_{2}(\Omega).
\end{eqnarray*}

\indent We need the following regularity assumption. For any $f\in
L_{2}(\Omega)$, $Tf\in H_{0}^{1}(\Omega)\cap
H^{1+\gamma_{1}}(\Omega)$ and $T^{*}f\in H_{0}^{1}(\Omega)\cap
H^{1+\gamma_{2}}(\Omega)$ satisfying
\begin{eqnarray}\label{s2.21}
&&\|Tf\|_{1+\gamma_{1}}\leq C_{\Omega}\|f\|_{0},\\\label{s2.22}
&&\|T^{*}f\|_{1+\gamma_{2}}\leq C_{\Omega}\|f\|_{0}.
\end{eqnarray}
According to
\cite{grisvard} and the section 5.5 in \cite{brenner}, the above assumption is reasonable.\\

\indent For some $G\subset\Omega$, we need the following local regularity assumption.\\

\indent{\bf R(G)}. For any $f\in L_{2}(G)$, there exists a $\phi\in
H_{0}^{1}(G)\cap H^{1+\gamma_1}(G)$ satisfying
$$ a(\phi,v)=b(f,v),~~~\forall v\in H_{0}^{1}(G),$$
and
\begin{eqnarray}\label{s2.23}
\|\phi\|_{1+\gamma_{1},G}\leq C_{G}\|f\|_{-1+\gamma_{1},G}.
\end{eqnarray}
For any $g\in L_{2}(G)$, there exists a $\varphi\in H_{0}^{1}(G)\cap
H^{1+\gamma_2}(G)$ satisfying
$$ a(v, \varphi)=b(v, g),~~~\forall v\in H_{0}^{1}(G),$$
and
\begin{eqnarray}\label{s2.24}
\|\varphi\|_{1+\gamma_{2},G}\leq C_{G}\|g\|_{-1+\gamma_{2},G}.
\end{eqnarray}
Where $C_{\Omega}$, $C_{G}$ are two priori constants, and not necessarily the same at different places.\\

\indent Define the Ritz projection $P_{h}: H_{0}^{1}(\Omega)\to
V_{h}^{0}(\Omega)$ and $P_{h}^{*}: H_{0}^{1}(\Omega)\to
V_{h}^{0}(\Omega)$ by
\begin{eqnarray}\label{s2.25}
a(u-P_{h}u, v)=0,~~and~~a(v, u-P_{h}^{*}u)=0,~~~\forall v\in
V_{h}^{0}(\Omega).
\end{eqnarray}
Then $T_{h}=P_{h}T$, $T_{h}^*=P_{h}^*T^*$ (see \cite{babuska}).

Let $M(\lambda)$ be the space spanned by all generalized
eigenfunctions corresponding to $\lambda$ of $T$, $M_{h}(\lambda)$
be the space spanned by all generalized eigenfunctions corresponding
to all eigenvalues of $T_{h}$ that converge to $\lambda$. In view of
the adjoint problem (\ref{s2.11}) and (\ref{s2.12}), the definitions
of $M^{*}(\lambda^{*})$ and $ M_{h}^{*}(\lambda^{*})$ are analogous
to $M(\lambda)$ and $M_{h}(\lambda)$.\\
\indent In this paper, we suppose that $\lambda$ is an eigenvalue of
(\ref{s2.3}) with the algebraic multiplicity $q$ and the ascent
$\alpha=1$. Then $\lambda^*=\overline{\lambda}$ be eigenvalue of
(\ref{s2.11}), $M(\lambda)$ and $M^{*}(\lambda^{*})$
are all eigenfunction space.\\
\indent Let $\lambda_{h}$ be the eigenvalue of (\ref{s2.8}) which
converges to $\lambda$, let
$\lambda_{h}^{*}=\overline{\lambda_{h}}$, and
$M^{*}(\lambda_{h}^{*})$ be the generalized eigenfunction space
corresponding to the eigenvalue $\lambda_{h}^{*}$ of $T_{h}^{*}$.\\

\indent {\bf Remark 2.1. }~~Obviously, it's difficult to determine the ascent $\alpha$ of the eigenvalue $\lambda$ of (\ref{s2.3}) theoretically. But one could easily find that when the ascents of the eigenvalues of (\ref{s2.8}), which converge to the same eigenvalue $\lambda$ of (\ref{s2.3}), are all equal to 1,
  one can conclude that the ascent $\alpha=1$  from the standard theory of spectral approximation. And the ascents of eigenvalues of (\ref{s2.8}) can be determined by computation.

We also need the lemma as follows (see \cite{kolman,yang2}):\\

\indent{\bf Lemma 2.1.}~~Let $(\lambda,u)$ be an eigenpair of
(\ref{s2.3}), and $(\lambda^{*}=\overline{\lambda},u^{*})$ be the
associated eigenpair of the adjoint problem (\ref{s2.11}). Then for
all $w,w^{*}\in H^{1}_{0}(\Omega)$, $b(w,w^{*})\neq 0$,
\begin{eqnarray}\label{s2.26}
\frac{a(w,w^{*})}{b(w,w^{*})}-\lambda=\frac{a(w-u,w^{*}-u^{*})}{b(w,w^{*})}-\lambda\frac{b(w-u,w^{*}-u^{*})}{b(w,w^{*})}.
\end{eqnarray}
\indent{\bf Proof.}~~see \cite{kolman,yang2}.~~~$\square$\\

\indent The a priori error estimates of the finite element
approximations (\ref{s2.8}) and (\ref{s2.12}) can refer to \cite{babuska,bramble}.\\

\indent{\bf Lemma 2.2.}~~Assume that $M(\lambda)\subset
H^{r+s}(\Omega)$, $M^*(\lambda^*)\subset H^{r+s_{2}}(\Omega)$
($0<s,s_{2}<1$). Then
\begin{eqnarray}\label{s2.27}
|\lambda_{h}-\lambda| \lesssim h^{r+s-1+r+s_{2}-1};
\end{eqnarray}
let $u_{h}\in M_{h}(\lambda)$ with $\|u_{h}\|_{0}=1$, then there is
$u\in M(\lambda)$ such that
\begin{eqnarray}\label{s2.28}
&&\|u_{h}-u\|_{1} \lesssim  h^{r+s-1},\\\label{s2.29}
&&\|u_{h}-u\|_{0} \lesssim  h^{r+s-1+\gamma_{2}};
\end{eqnarray}
let $u_{h}^{*}\in M_{h}^*(\lambda^*)$ with $\|u_{h}^*\|_{0}=1$, then
there is $u^{*}\in M^*(\lambda^*)$ such that
\begin{eqnarray}\label{s2.30}
&&\|u_{h}^{*}-u^{*}\|_{1} \lesssim h^{r+s_{2}-1};\\\label{s2.31}
&&\|u_{h}^{*}-u^{*}\|_{0} \lesssim h^{r+s_{2}-1+\gamma_{1}}.
\end{eqnarray}
\indent{\bf Proof.}~~see \cite{babuska}.~~~$\square$\\

\indent \cite{wahlbin,xu1} etc. studied the local behavior of finite
element. The following Lemma 2.3 is a simple generalization of Lemma
3.2 in \cite{xu1}. We can easily prove this Lemma by the same
argument as that of Lemma 3.2 in \cite{xu1}.\\

\indent {\bf Lemma 2.3.}~~Suppose that $f\in H^{-1}(\Omega)$ and
$G\subset\subset \Omega_{0}\subset \Omega$. If $w\in
V_{h}(\Omega_{0})$, $w|_{\partial\Omega\cap\partial \Omega_{0}}=0$,
satisfies
\begin{eqnarray}\label{s2.32yy}
a(w,v)=f(v),~~~\forall v\in V_{0}^{h}(\Omega_{0}),
\end{eqnarray}
then
\begin{eqnarray}\label{s2.32}
\|w\|_{1,G}\lesssim \|w\|_{0,\Omega_{0}}+\|f\|_{-1,\Omega_{0}},
\end{eqnarray}
where
\begin{eqnarray*}
\|f\|_{-1,\Omega_{0}}=\sup\limits_{\phi\in
H_{0}^{1}(\Omega_{0}),\|\phi\|_{1,\Omega_{0}}=1}f(\phi).
\end{eqnarray*}
\indent {\bf Proof.~~}
 Let $p\geq 2\nu-1$
be an integer, and let
 $$D\subset\subset \Omega_{p}\subset\subset \Omega_{p-1}\subset\subset\cdots\subset\subset \Omega_{1}\subset\subset \Omega_{0}.$$
 Choose $D_{1}\subset \Omega$ satisfying $D\subset\subset D_{1}\subset\subset \Omega_{p}$
and $\omega\in C^{\infty}(\bar{\Omega})$ such that
$\mathrm{supp}~\omega \subset\subset \Omega_{p}$ and $\omega\equiv
1$ on $\overline{D_{1}}$. Then, from (A3), there exists $v\in
V_{0}^{h}(\Omega_{p})$ such that
\begin{eqnarray*}
\|\omega^{2} w-v\|_{1,\Omega_{p}}\lesssim
h_{\Omega_{0}}\|w\|_{1,\Omega_{p}},
\end{eqnarray*}
so we have
\begin{eqnarray}\label{s3.5r}
a(w, \omega^{2} w-v)\lesssim h_{\Omega_{0}}\|w\|_{1,\Omega_{p}}^{2}
\end{eqnarray}
and
\begin{eqnarray}\label{s3.6r}
&&|f(v)|\lesssim
\|f\|_{-1,\Omega_{0}}\|v\|_{1,\Omega_{p}}\nonumber\\
&&~~~\lesssim
\|f\|_{-1,\Omega_{0}}(h_{\Omega_{0}}\|w\|_{1,\Omega_{p}}+\|\omega
w\|_{1,\Omega}).
\end{eqnarray}
Since $v\in V_{0}^{h}(\Omega_{p})\subset  V_{0}^{h}(\Omega_{0})$,
(\ref{s2.32yy}) implies
\begin{eqnarray}\label{s3.7r}
a(w, \omega^{2} w)=a(w, \omega^{2} w-v)+f(v).
\end{eqnarray}
Let
$a_{0}(u,v)=\int\limits_{\Omega}\sum\limits_{i,j=1}^{d}a_{ij}\partial_{i}u\partial_{j}\overline{v}$.
We can be derived from proof of Lemma 3.1 in \cite{xu1} that if
$\Omega_{0}\subset\Omega\subset R^{d}(d=2,3)$, $\omega\in
C^{\infty}(\bar{\Omega})$ , $\mathrm{supp}~\omega\subset\subset
\Omega_{0}$, then
\begin{eqnarray}\label{s3.8r}
a_{0}(\omega w,\omega w)\lesssim a(w,\omega^{2}
w)+\|w\|_{0,\Omega_{0}}^{2},~~~ \forall w\in H_{0}^{1}(\Omega).
\end{eqnarray}
It follows from (\ref{s3.5r})-(\ref{s3.8r}) that
\begin{eqnarray*}
&&\|\omega w\|_{1,\Omega}^{2}\lesssim a_{0}(\omega w,\omega w)\lesssim a(w,\omega^{2} w)+\|w\|_{0,\Omega_{0}}^{2}\\
&&=a(w,\omega^{2} w-v)+\|w\|_{0,\Omega_{0}}^{2}+f(v)\\
&&\lesssim
h_{\Omega_{0}}\|w\|_{1,\Omega_{p}}^{2}+\|w\|_{0,\Omega_{0}}^{2}+\|f\|_{-1,\Omega_{0}}(h_{\Omega_{0}}\|w\|_{1,\Omega_{p}}+\|\omega
w\|_{1,\Omega}),
\end{eqnarray*}
thus
\begin{eqnarray}\label{s3.9r}
\| w\|_{1,D}\lesssim
 h_{\Omega_{0}}^{1/2}\|w\|_{1,\Omega_{p}}+\|w\|_{0,\Omega_{0}}+\|f\|_{-1,\Omega_{0}}.
\end{eqnarray}
Similarly, we can get
\begin{eqnarray}\label{s3.10r}
\| w\|_{1,\Omega_{j}}\lesssim
 h_{\Omega_{0}}^{1/2}\|w\|_{1,\Omega_{j-1}}+\|w\|_{0,\Omega_{0}}+\|f\|_{-1,\Omega_{0}},~~~j=1,2,\cdots,p.
\end{eqnarray}
By using (\ref{s3.9r}) and (\ref{s3.10r}), we get from (A0) and (A2)
that
\begin{eqnarray*}
&&\| w\|_{1,D}\lesssim
 h_{\Omega_{0}}^{(p+1)/2}\|w\|_{1,\Omega_{0}}+\|w\|_{0,\Omega_{0}}+\|f\|_{-1,\Omega_{0}}\\
 &&\lesssim
 h_{\Omega_{0}}^{(p+1)/2}\|
h^{-1}w\|_{0,\Omega_{0}}+\|w\|_{0,\Omega_{0}}+\|f\|_{-1,\Omega_{0}}\\
  &&\lesssim
 \|w\|_{0,\Omega_{0}}+\|f\|_{-1,\Omega_{0}}.
\end{eqnarray*}
This completes the proof.
~~~$\square$\\

\indent {\bf Lemma 2.4.}~~Suppose that $G\subset\subset
\Omega_{0}\subset \Omega$. Then the following estimates are valid:
\begin{eqnarray}\label{s2.33}
&&h^{\gamma_{2}}\|u-P_{h}u\|_{1,\Omega}+\|u-P_{h}u\|_{0,\Omega}
\lesssim h^{\gamma_{2}}\inf\limits_{v\in V_{h}^{0}(\Omega)}
\|u-v\|_{1,\Omega},\\\label{s2.34}
 && \|u-P_{h}u\|_{1,G}\lesssim
\inf\limits_{v\in V_{h}^{0}(\Omega)}
\|u-v\|_{1,\Omega_{0}}+h^{\gamma_{2}}\|u-P_{h}u\|_{1,\Omega}.
\end{eqnarray}
\indent {\bf Proof.~~}For proof of (\ref{s2.33}) cf.
\cite{brenner,ciarlet}, for proof of (\ref{s2.34}) cf. Theorem 3.4
in \cite{xu1}.~~~$\square$\\

\indent {\bf Remark 2.2. }~~In \cite{xu1}, the condition {\em Superapproximation} is given as follows.\\
 \indent {\bf A.3.}~{\em Superapproximation.} For $G\subset
\Omega_{0}$, let $\omega\in C_{0}^{\infty}(\Omega)$ with
$\mathrm{supp}~\omega \subset\subset G$. Then for any $w\in
V_{h}(G)$, there exists $v\in V_{0}^{h}(G)$ such that
$$ \|h^{-1}(\omega w-v)\|_{1,G}\lesssim \|w\|_{1,G}.$$
In the proof of Lemma 3.2 in \cite{xu1} the authors choose $D_{1}\subset
\Omega$ satisfying $D\subset\subset D_{1}\subset\subset \Omega_{p}$
and $\omega\in C_{0}^{\infty}(\Omega)$ such that $\omega\equiv 1$ on
$\overline{D_{1}}$ and
 $\mathrm{supp}~\omega
\subset\subset \Omega_{p}$.\\
This paper just makes a minor modification, so that
the theory of the local error estimates built in \cite{xu1} applies to the local domains containing the corner points,
see Lemma 2.3 and Lemma 2.4.

\section{Multilevel discretizations based on local defect-correction}
Consider the eigenvalue problem (\ref{s2.3}) which has an isolated
singular point $z\in \overline{\Omega}$~(e.g., see ~Figure 3.1).\\
\indent Let $D\subset\subset \Omega$ be a given subdomain containing
the singular point $z$, and we introduce domains
$$\Omega\supset\Omega_{1}\supset\Omega_{2}\supset\cdots\supset\Omega_{l}\supset\supset D.$$
\indent Let $\pi_{H}(\Omega)$ be a shape-regular grid, which is made
up of simplices, with size $H\in(0,1)$, $\pi_{w}(\Omega)$ be a
refined mesoscopic shape-regular grid(from $\pi_{H}(\Omega)$) and
$\pi_{h}(\Omega_{i})$ be a locally refined grid (from
$\pi_{h_{i-1}}(\Omega_{i-1})$) that satisfy $h_{-1}=H$,
$h_{0}=w,~~~h_{i}\ll h_{i-1}~(i=0,1,\cdots,l)$. (Figure 3.1
 shows~$\pi_{H}(\Omega)$, $\pi_{w}(\Omega)$ and
$\pi_{h_{1}}(\Omega_{1})$). Let $V_{H}^{0}(\Omega),
V_{w}^{0}(\Omega)$, and $\{V_{h_{i}}^{0}(\Omega_{i})\}_{1}^{l}$ be
finite element spaces
 of degree less than or equal to $r$ defined on $\pi_{H}(\Omega)$, $\pi_{w}(\Omega)$
and $\{\pi_{h_{i}}(\Omega_{i})\}_{1}^{l}$), respectively.\\

\begin{center}
\includegraphics*[110,380][500,475]{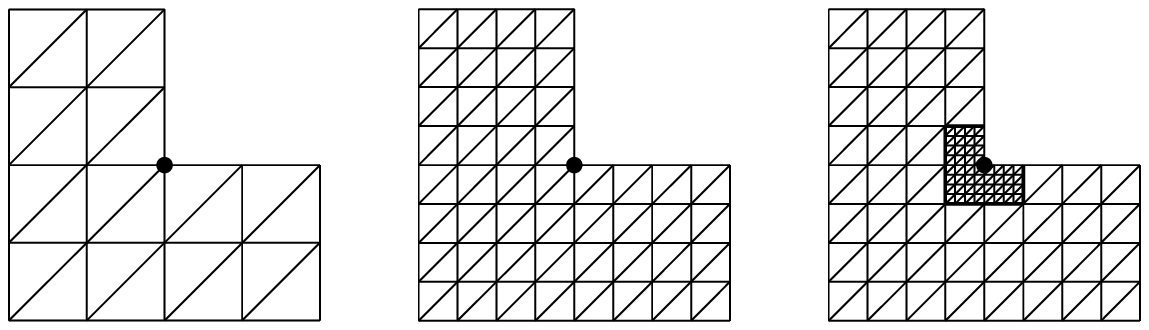}
\put(-220,-20){\bf \small Figure 3.1}
\end{center}

\indent Based on Algorithm $B_{0}$ in \cite{dai} we establish the
following
three-level discretization scheme.\\
\indent{\bf Scheme 3.1}(Three-level discretizations based on local defect-correction.).\\
\indent {\bf Step 1.} Solve (\ref{s2.3}) on a globally coarse grid
$\pi_{H}(\Omega)$: find $\lambda_{H}\in \mathcal{C}, u_{H}\in
V_{H}^{0}(\Omega)$ such that $\|u_{H}\|_{0}=1$ and
\begin{eqnarray*}
a(u_{H},v)=\lambda_{H} b(u_{H},v),~~~\forall v\in V_{H}^{0}(\Omega).
\end{eqnarray*}
Let~$\lambda_{H}^*=\overline{\lambda}_{H}$, and find $u_{H}^{*}\in
M^*(\lambda_{H}^*)$ with $\|u_{H}^{*}\|_{0}=1$ such that $|b(u_{H},
u_{H}^*)|$ has a positive lower bound uniformly with respect to $H$
(see section 5.1).\\
\indent {\bf Step 2.} Solve two linear boundary value problems on a
globally mesoscopic grid $\pi_{w}(\Omega)$: find $u^w\in
V_{w}^{0}(\Omega)$ such that
\begin{eqnarray*}
a(u^w, v)=\lambda_{H}b(u_{H},v),~~~\forall v\in V_{w}^{0}(\Omega);
\end{eqnarray*}
find $u^{w*}\in V_{w}^{0}(\Omega)$ such that
\begin{eqnarray*}
a(v, u^{w*})=\lambda_{H}b(v, u_{H}^{*}),~~~\forall v\in
V_{w}^{0}(\Omega).
\end{eqnarray*}
Then compute the Rayleigh
quotient $ \lambda^{w}=\frac{a(u^{w},u^{w*})}{ b(u^{w},u^{w*})}$.\\
\indent {\bf Step 3.} Solve two linear boundary value problems on a
locally fine grid $\pi_{h_{1}}(\Omega_{1})$: find $e^{h_{1}}\in
V_{h_{1}}^{0}(\Omega_{1})$ such that
\begin{eqnarray}\label{s3.1}
a(e^{h_{1}},v)= \lambda^{w}b(u^{w},v)-a(u^{w},v),~~~\forall v\in
V_{h_{1}}^{0}(\Omega_{1});
\end{eqnarray}
find $e^{h_{1}*}\in V_{h_{1}}^{0}(\Omega_{1})$ such that
\begin{eqnarray}\label{s3.2}
a(v, e^{h_{1}*})= \lambda^{w}b(v, u^{w*})-a(v, u^{w*}),~~~\forall
v\in V_{h_{1}}^{0}(\Omega_{1}).
\end{eqnarray}
\indent {\bf Step 4.} Set
\begin{eqnarray}\label{s3.3}
u^{w,h_{1}}= \left \{
\begin{array}{ll}
u^{w}+e^{h_{1}}&~~~on~ \overline{\Omega}_{1},\\
u^{w}&~~~in~ \Omega\setminus\overline{\Omega}_{1}
\end{array}
\right.
\end{eqnarray}
\begin{eqnarray}\label{s3.4}
u^{w,h_{1}*}= \left \{
\begin{array}{ll}
u^{w*}+e^{h_{1}*}&~~~on~ \overline{\Omega}_{1},\\
u^{w*}&~~~in~ \Omega\setminus\overline{\Omega}_{1}
\end{array}
\right.
\end{eqnarray}
And compute the Rayleigh quotient
\begin{eqnarray}\label{s3.5}
\lambda^{w,h_{1}}=\frac{a(u^{w,h_{1}},u^{w,h_{1}*})}{
b(u^{w,h_{1}},u^{w,h_{1}*})},~~~\lambda^{w,h_{1}*}=\overline{\lambda^{w,h_{1}}}.
\end{eqnarray}

\indent We use $(\lambda^{w,h_{1}}, u^{w,h_{1}})$ and
$(\lambda^{w,h_{1}*}, u^{w,h_{1}*})$ obtained by Scheme 3.1 as
the approximate eigenpair of (\ref{s2.3}) and (\ref{s2.11}), respectively.\\
\indent It is obvious that $(\lambda^{w},u^{w})$ and
$(\overline{\lambda^{w}}, u^{w*})$ in Scheme 3.1 can be viewed as
approximate eigenpairs obtained by the two-grid discretization
scheme in \cite{kolman,yang2} from $\pi_{H}(\Omega)$ and
$\pi_{w}(\Omega)$.\\

\indent  Using Scheme 3.1, abrupt changes of mesh size can  appear
near $\partial \Omega_{1}$. Influenced by the technique on the
transition layer proposed by \cite{bi}, we repeatedly use the local
defect-correction technique to establish the following
multilevel discretization scheme.\\
\indent{\bf Scheme 3.2}(multilevel discretizations based on local defect-correction.).\\
\indent {\bf Step 1.}The same as that of~Step 1 of Scheme 3.1.\\
\indent {\bf Step 2.}The same as that of~Step 2 of Scheme 3.1.\\
\indent {\bf Step 3.} $u^{w,h_{0}}\Leftarrow u^{w}$,
$\lambda^{w,h_{0}}\Leftarrow \lambda^{w}$,
 $u^{w,h_{0}*}\Leftarrow u^{w*}$, $\lambda^{w,h_{0}*}\Leftarrow \lambda^{w*}$.\\
\indent {\bf Step 4.} For $i=1,2,\cdots,l$, execute Step 5 and Step 6.\\
 \indent{\bf Step 5.} Solve linear
boundary value problems on locally fine grid
$\pi_{h_{i}}(\Omega_{i})$: find $e^{h_{i}}\in
V_{h_{i}}^{0}(\Omega_{i})$ such that
\begin{eqnarray}\label{s3.6}
a(e^{h_{i}},v)=
\lambda^{w,h_{i-1}}b(u^{w,h_{i-1}},v)-a(u^{w,h_{i-1}},v),~~~\forall
v\in V_{h_{i}}^{0}(\Omega_{i});
\end{eqnarray}
find $e^{h_{i}*}\in V_{h_{i}}^{0}(\Omega_{i})$ such that
\begin{eqnarray}\label{s3.7}
a(v, e^{h_{i}*})= \lambda^{w,h_{i-1}*}b(v, u^{w,h_{i-1}*})-a(v,
u^{w,h_{i-1}*}),~~~\forall v\in V_{h_{i}}^{0}(\Omega_{i}).
\end{eqnarray}
\indent {\bf Step 6.} Set
\begin{eqnarray}\label{s3.8}
u^{w,h_{i}}= \left \{
\begin{array}{ll}
u^{w,h_{i-1}}+e^{h_{i}}&~~~on~ \overline{\Omega}_{i},\\
u^{w,h_{i-1}}&~~~in~ \Omega\setminus\overline{\Omega}_{i}
\end{array}
\right.
\end{eqnarray}
\begin{eqnarray}\label{s3.9}
u^{w,h_{i}*}= \left \{
\begin{array}{ll}
u^{w,h_{i-1}*}+e^{h_{i}*}&~~~on~ \overline{\Omega}_{i},\\
u^{w,h_{i-1}*}&~~~in~ \Omega\setminus\overline{\Omega}_{i}
\end{array}
\right.
\end{eqnarray}
And compute
\begin{eqnarray}\label{s3.10}
 \lambda^{w,h_{i}}=\frac{a(u^{w,h_{i}},u^{w,h_{i}*})}{
b(u^{w,h_{i}},u^{w,h_{i}}*)},~~~\lambda^{w,h_{i}*}=\overline{\lambda^{w,h_{i}}}.
\end{eqnarray}

\indent We use $(\lambda^{w,h_{l}}, u^{w,h_{l}})$ and
$(\lambda^{w,h_{l}*}, u^{w,h_{l}*})$ obtained by Scheme 3.2 as the
approximate eigenpair of (\ref{s2.3}) and (\ref{s2.11}),
respectively.

\section{Theoretical Analysis}

\indent Next we shall discuss the error estimates of Scheme 3.1 and Scheme 3.2.\\
\indent In our analysis, we introduce an auxiliary grid
$\pi_{h_{i}}(\Omega)$ which is defined globally, and denote the
piecewise polynomials space of degree $\leq r$ by
$V_{h_{i}}^{0}(\Omega)$ ($i=1,2,\cdots,l$). We also assume that
$\pi_{h_{i}}(\Omega_{i})$ and $V_{h_{i}}^{0}(\Omega_{i})$ are the
restrictions of $\pi_{h_{i}}(\Omega)$ and $V_{h_{i}}^{0}(\Omega)$ to
$\Omega_{i}$, respectively, and
$$ V_{H}^{0}(\Omega)\subset V_{w}^{0}(\Omega)\subset
V_{h_{1}}^{0}(\Omega)\subset
V_{h_{2}}^{0}(\Omega)\subset\cdots\subset  V_{h_{l}}^{0}(\Omega).$$

\indent For $D$ and $\Omega_{i}$ stated at the beginning of section
3, let $G_{i}\subset \Omega$ and $F\subset \Omega$ satisfy
$D\subset\subset F\subset\subset G_{i}\subset\subset\Omega_{i}$
($i=1,2,\cdots,l$).\\

\indent{\bf Theorem 4.1.}~~Assume that $M(\lambda)\subset
H_{0}^{1}(\Omega)\cap H^{r+s}(\Omega)\cap
H^{r+1}(\Omega\setminus\overline{D})$ and ($1<r+s$, $0\leq s<1$),
 $M^*(\lambda^*)\subset
H_{0}^{1}(\Omega)\cap H^{r+s_{2}}(\Omega)\cap
H^{r+1}(\Omega\setminus\overline{D})$ and ($1<r+s_{2}$, $0\leq
s_{2}<1$), and $H$ is properly small. Then there exists $u\in
M(\lambda)$ and $u^{*}\in M^*(\lambda^*)$ such that
\begin{eqnarray}\label{s4.1}
\|u^{w}-u\|_{1}&\lesssim &
H^{r+s-1+\gamma_{2}}+w^{r+s-1},\\\label{s4.2}
 \|u^{w}-u\|_{0}&\lesssim &
H^{r+s-1+\gamma_{2}},\\\label{s4.3}
 \|u^{w}-u\|_{1,
\Omega\setminus\overline{F}}&\lesssim &
H^{r+s-1+\gamma_{2}}+w^{r},\\\label{s4.4}
\|u^{w*}-u^{*}\|_{1}&\lesssim &
H^{r+s_{2}-1+\gamma_{1}}+w^{r+s_{2}-1},\\\label{s4.5}
 \|u^{w*}-u^{*}\|_{0}&\lesssim &
H^{r+s_{2}-1+\gamma_{1}},\\\label{s4.6}
\|u^{w*}-u^{*}\|_{1,\Omega\setminus\overline{F}}&\lesssim &
H^{r+s_{2}-1+\gamma_{1}}+w^{r},\\\label{s4.7}
 \mid \lambda^{w}-\lambda
\mid&\lesssim &
H^{2r+s+s_{2}-2+\gamma_{1}+\gamma_{2}}+w^{2r+s+s_{2}-2}.
\end{eqnarray}
\indent {\bf Proof.}~~Let~$u\in M(\lambda)$ and $u^*\in
M^*(\lambda^*)$  such that~$u-u_{H}$ and $u^*-u_{H}^*$ both satisfy Lemma~2.2. From
(\ref{s2.13}), (\ref{s2.14}), Step 2 of Scheme 3.1, (\ref{s2.15}),
Lemma 2.2 and Lemma 2.4, we derive that
\begin{eqnarray*}
&&\|u^{w}-u\|_{1}=\|\lambda_{H}T_{w}u_{H}-\lambda T u\|_{1}\nonumber\\
&&~~~\leq\|\lambda_{H}T_{w}u_{H}-\lambda T_{w} u\|_{1}+\|\lambda
T_{w} u-\lambda T u\|_{1}\nonumber\\
&&~~~\lesssim  \|\lambda_{H}u_{H}-\lambda  u\|_{0}+\lambda \|
P_{w}Tu-T u\|_{1}\nonumber\\
&&~~~\lesssim  H^{r+s-1+\gamma_{2}}+w^{r+s-1},
\end{eqnarray*}
then (\ref{s4.1}) follows. By Lemma 2.2 and Lemma 2.4,
\begin{eqnarray*}
&&\|u^{w}-u\|_{1, \Omega\setminus\overline{F}} \lesssim
\|\lambda_{H}u_{H}-\lambda u\|_{0}+\lambda \|
P_{w}Tu-T u\|_{1, \Omega\setminus\overline{F}}\nonumber\\
&&~~~\lesssim  H^{r+s-1+\gamma_{2}}+w^{r},
\end{eqnarray*}
then (\ref{s4.3}) follows. By calculation,
\begin{eqnarray*}
&&\|u^{w}-u\|_{0}=\|\lambda_{H}T_{w}u_{H}-\lambda T u\|_{0}\nonumber\\
&&~~~\leq\|\lambda_{H}T_{w}u_{H}-\lambda T_{w} u\|_{0}+\|\lambda
T_{w} u-\lambda T u\|_{0}\nonumber\\
&&~~~\lesssim  \|\lambda_{H}u_{H}-\lambda  u\|_{0}+\lambda \|
P_{w}Tu-T u\|_{0}\nonumber\\
&&~~~\lesssim  H^{r+s-1+\gamma_{2}}+w^{r+s-1+\gamma_{2}}\nonumber\\
&&~~~\lesssim  H^{r+s-1+\gamma_{2}},
\end{eqnarray*}
then (\ref{s4.2}) follows.\\
Similarly we can prove~(\ref{s4.4}), (\ref{s4.5}) and (\ref{s4.6}). From
(\ref{s2.26}), we have
\begin{eqnarray}\label{s4.8}
\lambda^{w}-\lambda=\frac{a(u^{w}-u, u^{w*}-u^{*})}{b(u^{w},
u^{w*})}-\lambda\frac{b(u^{w}-u, u^{w*}-u^{*})}{b(u^{w}, u^{w*})}.
\end{eqnarray}
Note that~$u_{H}$  and~$u^{w}$ just approximate the same eigenfuntion $u$, $u_{H}^{*}$ and~$u^{w*}$
approximate the same adjoint eigenfuntion $u^{*}$, $|b(u_{H}, u_{H}^*)|$ has a positive lower
bound uniformly with respect to $H$, therefore~$b(u^{w}, u^{w*})$ has a
positive lower bound uniformly. Combining ~(\ref{s4.1}), (\ref{s4.2}),
(\ref{s4.4}), (\ref{s4.5}) and~(\ref{s4.8}) yields~(\ref{s4.7}).~~~$\square$\\

The following Theorem 4.2 is a critical result in this paper, which develops the results of Theorem 3.3 in \cite{dai}.\\

\indent {\bf Theorem 4.2.}~~ Assume that
$R(\Omega_{i})$ holds ($i=1,2,\cdots,l$), $u\in M(\lambda)$ and
$u^{*}\in M^*(\lambda^*)$. Then
\begin{eqnarray}\label{s4.9}
&&\|u^{w,h_{l}}-P_{h_{l}}u\|_{1,\Omega}\lesssim
 \|u-P_{h_{l}}u\|_{0,\Omega_{l}}+
h_{l-1}^{\gamma_{2}}\|P_{h_{l}}u-u^{w,h_{l-1}}\|_{1,\Omega_{l}}\nonumber\\
&&~~~~~~~+\|\lambda u -\lambda^{w,h_{l-2}} u^{w,h_{l-2}}\|_{0,
\Omega_{l}}+\|\lambda^{w,h_{l-1}}u^{w,h_{l-1}}-\lambda
u\|_{0}\nonumber\\
&&~~~~~~+\|u^{w,h_{l-1}}-P_{h_{l}}u\|_{1,\Omega\backslash\overline{G_{l}}}+
\|u^{w,h_{l-1}}-u\|_{1, \Omega_{l}\backslash \overline{F}},~~~l\geq
1;
\end{eqnarray}
\begin{eqnarray}\label{s4.10}
&&\|u^{w,h_{l}*}-P_{h_{l}}^*u\|_{1,\Omega}\lesssim
 \|u^*-P_{h_{l}}^*u^*\|_{0,\Omega_{l}}+
h_{l-1}^{\gamma_{1}}\|P_{h_{l}}^*u^*-u^{w,h_{l-1}*}\|_{1,\Omega_{l}}\nonumber\\
&&~~~~~~~+\|\lambda u^* -\lambda^{w,h_{l-2}*} u^{w,h_{l-2}*}\|_{0,
\Omega_{l}}+\|\lambda^{w,h_{l-1}*}u^{w,h_{l-1}*}-\lambda
u^*\|_{0}\nonumber\\
&&~~~~~~+\|u^{w,h_{l-1}*}-P_{h_{l}}u^*\|_{1,\Omega\backslash\overline{G_{l}}}+
\|u^{w,h_{l-1}*}-u^*\|_{1, \Omega_{l}\backslash
\overline{F}},~~~l\geq 1.
\end{eqnarray}
\indent {\bf Proof.}~~Due to the inequality
\begin{eqnarray}\label{s4.11}
&&\|u^{w,h_{l}}-P_{h_{l}}u\|_{1,\Omega} \lesssim
\|u^{w,h_{l}}-P_{h_{l}}u\|_{1,D}+\|u^{w,h_{l}}-P_{h_{l}}u\|_{1,G_{l}\backslash
\overline{D}}
\nonumber\\
&&~~~~~~+\|u^{w,h_{l}}-P_{h_{l}}u\|_{1,\Omega\backslash\overline{G_{l}}},
\end{eqnarray}
we shall estimate $\|u^{w,h_{l}}-P_{h_{l}}u\|_{1,D}$,
$\|u^{w,h_{l}}-P_{h_{l}}u\|_{1,G_{l}\backslash \overline{D}}$, and
$\|u^{w,h_{l}}-P_{h_{l}}u\|_{1,\Omega\backslash\overline{G_{l}}}$
, respectively.\\
\indent First, we proceed to estimate
$\|u^{w,h_{l}}-P_{h_{l}}u\|_{1,D}$. From (\ref{s3.8}), (\ref{s3.6})
and (\ref{s2.25}) we derive
\begin{eqnarray}\label{s4.12}
&&a(u^{w,h_{l}}-P_{h_{l}}u,v)
=a(u^{w,h_{l}},v)-a(P_{h_{l}}u,v)=a(u^{w,h_{l-1}}+e^{h_{l}},v)-a(u,v)\nonumber\\
&&~~~=\lambda^{w,h_{l-1}}b(u^{w,h_{l-1}},v) -\lambda b(u,v),~~~
\forall v \in V_{h_{l}}^{0}(\Omega_{l}).
\end{eqnarray}
It is obvious that
\begin{eqnarray}\label{s4.13}
&&~~~\lambda^{w,h_{l-1}}b(u^{w,h_{l-1}},v) -\lambda
b(u,v)\nonumber\\
&&=(\lambda^{w,h_{l-1}}-\lambda)b(u,v) +\lambda^{w,h_{l-1}}
b(u^{w,h_{l-1}}-u,v),~~~ \forall v \in H_{0}^{1}(\Omega),
\end{eqnarray}
which together with (\ref{s4.12}) yields
\begin{eqnarray*}
a(u^{w,h_{l}}-P_{h_{l}}u,v)=(\lambda^{w,h_{l-1}}-\lambda)b(u,v)
+\lambda^{w,h_{l-1}} b(u^{w,h_{l-1}}-u,v),~~~ \forall v \in
V_{h_{l}}^{0}(\Omega_{l}).
\end{eqnarray*}
Since
$V_{0}^{h_{l}}(\Omega_{l})\subset V_{h_{l}}^{0}(\Omega_{l})$ , thus,
from the above formula ~and Lemma 2.3 we deduce that
\begin{eqnarray}\label{s4.14}
\|u^{w,h_{l}}-P_{h_{l}}u\|_{1,D}\lesssim
\|u^{w,h_{l}}-P_{h_{l}}u\|_{0,\Omega_{l}}+
|\lambda^{w,h_{l-1}}-\lambda|+ \|u^{w,h_{l}}-u\|_{0,\Omega_{l}}.
\end{eqnarray}
By calculation, we have
\begin{eqnarray*}
&&\|u^{w,h_{l}}-P_{h_{l}}u\|_{0,\Omega_{l}}\leq
\|u^{w,h_{l-1}}-P_{h_{l}}u\|_{0,\Omega_{l}}+\|e^{h_{l}}\|_{0,\Omega_{l}}\\
&&~~~\leq
\|u-P_{h_{l}}u\|_{0,\Omega_{l}}+\|u-u^{w,h_{l-1}}\|_{0,\Omega_{l}}+\|e^{h_{l}}\|_{0,\Omega_{l}},
\end{eqnarray*}
substituting the above relation into (\ref{s4.14}) we obtain
\begin{eqnarray}\label{s4.15}
&&~~~\|u^{w,h_{l}}-P_{h_{l}}u\|_{1,D}\nonumber\\
&&\lesssim |\lambda^{w,h_{l-1}}-\lambda|+
\|u^{w,h_{l-1}}-u\|_{0,\Omega_{l}}+\|u-P_{h_{l}}u\|_{0,\Omega_{l}}+\|e^{h_{l}}\|_{0,\Omega_{l}}.
\end{eqnarray}
To estimate $\|e^{h_{l}}\|_{0,\Omega_{l}}$, we use the Aubin-Nitsche
duality argument. For any given $f\in L_{2}(\Omega_{l})$, consider
the boundary value problem: find $\varphi\in H_{0}^{1}(\Omega_{l})$
such that
\begin{eqnarray}\label{s4.16}
a(v, \varphi)=b(v, f)~~~\forall v\in H_{0}^{1}(\Omega_{l}).
\end{eqnarray}
Let $\varphi$ be the generalized solution of (\ref{s4.16}),
$\varphi_{h_{l}}$ and $\varphi_{h_{l-1}}$ be finite element
solutions of (\ref{s4.16}) in $V_{h_{l}}^{0}(\Omega_{l})$ and
$V_{h_{l-1}}^{0}(\Omega_{l})$, respectively. Then,
\begin{eqnarray}\label{s4.17}
\|\varphi-\varphi_{h_{l}}\|_{1,\Omega_{l}}\lesssim
h_{l}^{\gamma_{2}}\|f\|_{0,\Omega_{l}},~~~
\|\varphi-\varphi_{h_{l-1}}\|_{1,\Omega_{l}}\lesssim
h_{l-1}^{\gamma_{2}}\|f\|_{0,\Omega_{l}}.
\end{eqnarray}
From (\ref{s3.6}) and (\ref{s3.8}) we get
\begin{eqnarray*}
a(u^{w,h_{l}}, \varphi_{h_{l}})=\lambda^{w,h_{l-1}}b(u^{w,h_{l-1}},
\varphi_{h_{l}}),
\end{eqnarray*}
thus by the definitions of $\varphi$, $\varphi_{h_{l}}$ and
$e^{h_{l}}$, we deduce that
\begin{eqnarray}\label{s4.18}
&&b(e^{h_{l}}, f)=a(e^{h_{l}}, \varphi)=a(e^{h_{l}}, \varphi_{h_{l}})=a(u^{w,h_{l}}-u^{w,h_{l-1}}, \varphi_{h_{l}})\nonumber\\
&&~~~=a(P_{h_{l}}u-u^{w,h_{l-1}}, \varphi_{h_{l}})+a(u^{w,h_{l}}, \varphi_{h_{l}})-a(P_{h_{l}}u, \varphi_{h_{l}})\nonumber\\
&&~~~=a(P_{h_{l}}u-u^{w,h_{l-1}}, \varphi_{h_{l}})+\lambda^{w,h_{l-1}}b(u^{w,h_{l-1}}, \varphi_{h_{l}})-\lambda b(u, \varphi_{h_{l}})\nonumber\\
&&~~~=a(P_{h_{l}}u-u^{w,h_{l-1}}, \varphi_{h_{l}}-\varphi)+a(P_{h_{l}}u-u^{w,h_{l-1}}, \varphi-\varphi_{h_{l-1}})\nonumber\\
&&~~~+a(P_{h_{l}}u-u^{w,h_{l-1}}, \varphi_{h_{l-1}})
+\lambda^{w,h_{l-1}}b(u^{w,h_{l-1}}, \varphi_{h_{l}})-\lambda b(u, \varphi_{h_{l}})\nonumber\\
&&~~~\lesssim
h_{l-1}^{\gamma_{2}}\|P_{h_{l}}u-u^{w,h_{l-1}}\|_{1,\Omega_{l}}\|f\|_{0,
\Omega_{l}}+a(P_{h_{l}}u-u^{w,h_{l-1}},
\varphi_{h_{l-1}})\nonumber\\
&&~~~~~~+\lambda^{w,h_{l-1}}b(u^{w,h_{l-1}},
\varphi_{h_{l}})-\lambda b(u, \varphi_{h_{l}}).
\end{eqnarray}
Step 2 of Scheme 3.2 shows that
\begin{eqnarray*}
a(u^{w,h_{0}}, \varphi_{h_{0}})=\lambda^{w,h_{-1}}b(u^{w,h_{-1}},
\varphi_{h_{0}}),
\end{eqnarray*}
namely, for $l=1$,
\begin{eqnarray*}
a(u^{w,h_{l-1}},
\varphi_{h_{l-1}})=\lambda^{w,h_{l-2}}b(u^{w,h_{l-2}},
\varphi_{h_{l-1}}),
\end{eqnarray*}
for $l>1$, the above formula follows from (\ref{s3.6}) and (\ref{s3.8}). Therefore,
\begin{eqnarray*}
&&a(P_{h_{l}}u-u^{w,h_{l-1}}, \varphi_{h_{l-1}})=a(u-u^{w,h_{l-1}},
\varphi_{h_{l-1}})\nonumber\\
&&~~~=\lambda b(u, \varphi_{h_{l-1}})-a(u^{w,h_{l-1}}, \varphi_{h_{l-1}})\nonumber\\
&&~~~=\lambda b(u, \varphi_{h_{l-1}})-\lambda^{w,h_{l-2}} b(u^{w,h_{l-2}}, \varphi_{h_{l-1}})\nonumber\\
&&~~~\lesssim \|\lambda u -\lambda^{w,h_{l-2}} u^{w,h_{l-2}}\|_{0,
\Omega_{l}}\|f\|_{0,\Omega_{l}}.
\end{eqnarray*}
It is clear that
\begin{eqnarray*}
|\lambda^{w,h_{l-1}}b(u^{w,h_{l-1}}, \varphi_{h_{l}})-\lambda b(u,
\varphi_{h_{l}})|\lesssim \|\lambda^{w,h_{l-1}}u^{w,h_{l-1}}-\lambda
u\|_{0, \Omega_{l}}\|f\|_{0, \Omega_{l}}.
\end{eqnarray*}
Substituting the above two formulae into (\ref{s4.18}), we derive
\begin{eqnarray*}
|b(e^{h_{l}}, f)| &\lesssim &
(h_{l-1}^{\gamma_{2}}\|P_{h_{l}}u-u^{w,h_{l-1}}\|_{1,\Omega_{l}}+\|\lambda
u -\lambda^{w,h_{l-2}} u^{w,h_{l-2}}\|_{0,
\Omega_{l}}\nonumber\\
&&~~~+\|\lambda^{w,h_{l-1}}u^{w,h_{l-1}}-\lambda u\|_{0})\|f\|_{0,
\Omega_{l}}.
\end{eqnarray*}
Thus, we get
\begin{eqnarray}\label{s4.19}
\|e^{h_{l}}\|_{0,\Omega_{l}} &\lesssim &
h_{l-1}^{\gamma_{2}}\|P_{h_{l}}u-u^{w,h_{l-1}}\|_{1,\Omega_{l}}+\|\lambda
u -\lambda^{w,h_{l-2}} u^{w,h_{l-2}}\|_{0,
\Omega_{l}}\nonumber\\
&&~~~+\|\lambda^{w,h_{l-1}}u^{w,h_{l-1}}-\lambda u\|_{0}.
\end{eqnarray}
Substituting (\ref{s4.19}) into (\ref{s4.15}), we obtain
\begin{eqnarray}\label{s4.20}
&&\|u^{w,h_{l}}-P_{h_{l}}u\|_{1,D}\lesssim
\|u-P_{h_{l}}u\|_{0,\Omega_{l}}+
h_{l-1}^{\gamma_{2}}\|P_{h_{l}}u-u^{w,h_{l-1}}\|_{1,\Omega_{l}}\nonumber\\
&&~~~~~~~+\|\lambda u -\lambda^{w,h_{l-2}} u^{w,h_{l-2}}\|_{0,
\Omega_{l}}+|\lambda^{w,h_{l-1}}-\lambda |+\|u^{w,h_{l-1}}-
u\|_{0}.~~~~~
\end{eqnarray}
\indent Similarly, since $(G_{l}\setminus\overline{D})
\subset\subset \Omega_{l}$, we deduce
\begin{eqnarray}\label{s4.21}
&&\|u^{w,h_{l}}-P_{h_{l}}u\|_{1,D}\lesssim
\|u-P_{h_{l}}u\|_{0,\Omega_{l}}+
h_{l-1}^{\gamma_{2}}\|P_{h_{l}}u-u^{w,h_{l-1}}\|_{1,\Omega_{l}}\nonumber\\
&&~~~~~~~+\|\lambda u -\lambda^{w,h_{l-2}} u^{w,h_{l-2}}\|_{0,
\Omega_{l}}+|\lambda^{w,h_{l-1}}-\lambda |+\|u^{w,h_{l-1}}-
u\|_{0}.~~~~~
\end{eqnarray}
\indent The remainder is to analyze
$\|u^{w,h_{l}}-P_{h_{l}}u\|_{1,\Omega\backslash\overline{G}}$. From
(\ref{s3.8}), we see that
\begin{eqnarray*}
\|u^{w,h_{l}}-P_{h_{l}}u\|_{1,\Omega\backslash\overline{\Omega}_{l}}=\|u^{w,h_{l-1}}-P_{h_{l}}u\|_{1,\Omega\backslash\overline{\Omega}_{l}},
\end{eqnarray*}
which leads to
\begin{eqnarray}\label{s4.22}
&&~~~\|u^{w,h_{l}}-P_{h_{l}}u\|_{1,\Omega\backslash\overline{G_{l}}}\nonumber\\
&&\leq\|u^{w,h_{l}}-P_{h_{l}}u\|_{1,\Omega\backslash\overline{\Omega}_{l}}
+\|u^{w,h_{l-1}}-P_{h_{l}}u\|_{1,\Omega_{l}\backslash\overline{G_{l}}}+\|e^{h_{l}}\|_{1,\Omega_{l}\backslash\overline{G_{l}}}\nonumber\\
&&\lesssim
\|u^{w,h_{l-1}}-P_{h_{l}}u\|_{1,\Omega\backslash\overline{G_{l}}}+\|e^{h_{l}}\|_{1,\Omega_{l}\backslash\overline{G_{l}}}.
\end{eqnarray}
It follows from (\ref{s3.6}), (\ref{s2.3}) and (\ref{s4.13}) that
\begin{eqnarray*}
&&a(e^{h_{l}},v)=
\lambda^{w,h_{l-1}}b(u^{w,h_{l-1}},v)-a(u^{w,h_{l-1}},v) -\lambda
b(u,v)+a(u,v)\nonumber\\
&&~~~=(\lambda^{w,h_{l-1}}-\lambda)b(u,v) +\lambda^{w,h_{l-1}}
b(u^{w,h_{l-1}}-u,v)\nonumber\\
&&~~~~~~-a(u^{w,h_{l-1}}-u,v),~~\forall v\in V_{h}^{0}(\Omega_{l}),
\end{eqnarray*}
then, by Lemma 2.3, we have
\begin{eqnarray}\label{s4.23}
\|e^{h_{l}}\|_{1,\Omega_{l}\backslash \overline{G_{l}}}\lesssim
\|e^{h_{l}}\|_{0,\Omega_{l}\backslash
\overline{F_{}}}+|\lambda^{w,h_{l-1}}-\lambda|
+\|u^{w,h_{l-1}}-u\|_{1, \Omega_{l}\backslash \overline{F}},
\end{eqnarray}
where $F\subset\Omega$ satisfies $D\subset\subset F\subset\subset
G_{l}$. Substituting (\ref{s4.23}) into (\ref{s4.22}) we get
\begin{eqnarray*}
&&~~~\|u^{w,h_{l}}-P_{h_{l}}u\|_{1,\Omega\backslash\overline{G_{l}}}\nonumber\\
&&\lesssim
\|u^{w,h_{l-1}}-P_{h_{l}}u\|_{1,\Omega\backslash\overline{G_{l}}}+\|e^{h_{l}}\|_{0,\Omega_{l}\backslash
\overline{F}}+|\lambda^{w,h_{l-1}}-\lambda| +\|u^{w,h_{l-1}}-u\|_{1,
\Omega_{l}\backslash \overline{F}}.
\end{eqnarray*}
It follows from substituting (\ref{s4.19}) into the above inequality
that
\begin{eqnarray}\label{s4.24}
&&\|u^{w,h_{l}}-P_{h_{l}}u\|_{1,\Omega\backslash\overline{G_{l}}}
\lesssim
\|u^{w,h_{l-1}}-P_{h_{l}}u\|_{1,\Omega\backslash\overline{G_{l}}} +
h_{l-1}^{\gamma_{2}}\|P_{h_{l}}u-u^{w,h_{l-1}}\|_{1,\Omega_{l}}\nonumber\\
&&~~~~~~+\|\lambda u -\lambda^{w,h_{l-2}} u^{w,h_{l-2}}\|_{0,
\Omega_{l}}+\|\lambda^{w,h_{l-1}}u^{w,h_{l-1}}-\lambda
u\|_{0}\nonumber\\
&&~~~~~~+|\lambda^{w,h_{l-1}}-\lambda| +\|u^{w,h_{l-1}}-u\|_{1,
\Omega_{l}\backslash \overline{F}}.
\end{eqnarray}
\indent  Combining (\ref{s4.24}), (\ref{s4.20}), (\ref{s4.21}) and
(\ref{s4.11}),
finally, we obtain (\ref{s4.9}).\\
\indent We can prove~(\ref{s4.10}) by using the similar argument.
~~~$\square$\\

\indent {\bf Theorem 4.3.}~~Assume that the conditions of Theorem 4.1 hold. Then
there exists $u\in M(\lambda)$ and $u^{*}\in M^*(\lambda^*)$ such that
\begin{eqnarray}\label{s4.25}
&&\|u^{w,h_{1}}-u\|_{1,\Omega}\lesssim
h_{1}^{r+s-1}+w^{r}+H^{r+s-1+\gamma_{2}},\\\label{s4.26}
&&\|u^{w,h_{1}}-u\|_{0,\Omega}\lesssim
w^{r}+H^{r+s-1+\gamma_{2}},\\\label{s4.27}
&&\|u^{w,h_{1}}-u\|_{1,\Omega\setminus\overline{F}}\lesssim
w^{r}+H^{r+s-1+\gamma_{2}},\\\label{s4.28}
&&\|u^{w,h_{1}*}-u^*\|_{1,\Omega}\lesssim
h_{1}^{r+s_{2}-1}+w^{r}+H^{r+s_{2}-1+\gamma_{1}},\\\label{s4.29}
&&\|u^{w,h_{1}*}-u^*\|_{0,\Omega}\lesssim
w^{r}+H^{r+s_{2}-1+\gamma_{1}},\\\label{s4.30}
&&\|u^{w,h_{1}*}-u^*\|_{1,\Omega\setminus\overline{F}}\lesssim
w^{r}+H^{r+s_{2}-1+\gamma_{1}},\\\label{s4.31}
&&|\lambda^{w,h_{1}}-\lambda|\lesssim
h_1^{2r+s+s_{2}-2}+w^{2r}+H^{2r+s+s_{2}-2+\gamma_{1}+\gamma_{2}}.
\end{eqnarray}
\indent {\bf Proof.}~~Let $u\in M(\lambda)$ and $u^*\in
M^*(\lambda^*)$ such that~$u-u_{H}$ and $u^*-u_{H}^*$ both satisfy Lemma 2.2.
In Theorem 4.2,  choose $l=1$, $h_{-1}=H, h_{0}=w$, $u^{w,h_{0}}=u^{w}$,
$\lambda^{w,h_{0}}=\lambda^{w}$, $u^{w,h_{-1}}=u_{H}$,
$\lambda^{w,h_{-1}}=\lambda_{H}$, then we get
\begin{eqnarray}\label{s4.32}
&&\|u^{w,h_{1}}-P_{h_{1}}u\|_{1,\Omega}\lesssim
 \|u-P_{h_{1}}u\|_{0,\Omega_{1}}+
w^{\gamma_{2}}\|P_{h_{1}}u-u^{w}\|_{1,\Omega_{1}}\nonumber\\
&&~~~~~~~+\|\lambda u -\lambda_{H} u_{H}\|_{0,
\Omega_{1}}+\|\lambda^{w}u^{w}-\lambda
u\|_{0}\nonumber\\
&&~~~~~~+\|u^{w}-P_{h_{1}}u\|_{1,\Omega\backslash\overline{G_{1}}}+
\|u^{w}-u\|_{1, \Omega_{1}\backslash \overline{F_{}}}.
\end{eqnarray}
Using~Lemma 2.4, Theorem 4.1, Lemma 2.2 to estimate the terms at the right hand side of the above formula gives
\begin{eqnarray}\label{s4.33}
&&\|u^{w,h_{1}}-P_{h_{1}}u\|_{1,\Omega}\lesssim
h_{1}^{r+s-1+\gamma_{2}}+w^{\gamma_{2}}w^{r+s-1}+H^{r+s-1+\gamma_{2}}+w^{r+s-1+\gamma_{2}}\nonumber\\
&&~~~~~~+(w^{r+s-1+\gamma_{2}}+w^{r})+(w^{r+s-1+\gamma_{2}}+w^{r})\lesssim
H^{r+s-1+\gamma_{2}}+w^{r}.
\end{eqnarray}
Combining~(\ref{s2.33}) and (\ref{s2.34}) yields ~(\ref{s4.25}), (\ref{s4.26})
and (\ref{s4.27}). By the same argument we can prove~(\ref{s4.28}), (\ref{s4.29}) and
(\ref{s4.30}). From (\ref{s2.26}), we have
\begin{eqnarray}\label{s4.34}
\lambda^{w,h_{1}}-\lambda=\frac{a(u^{w,h_{1}}-u,
u^{w,h_{1}*}-u^{*})}{b(u^{w,h_{1}},
u^{w,h_{1}*})}-\lambda\frac{b(u^{w,h_{1}}-u,
u^{w,h_{1}*}-u^{*})}{b(u^{w,h_{1}}, u^{w,h_{1}*})}.
\end{eqnarray}
Note that $u_{H}$ and~$u^{w,h_{1}}$ just approximate the same
eigenfunction~$u$, $u_{H}^{*}$ and~$u^{w,h_{1}*}$ approximate the
same adjoint eigenfunction~$u^{*}$, $|b(u_{H}, u_{H}^*)|$ has a
positive lower bound uniformly with respect to $H$,
thus~$b(u^{w,h_{1}}, u^{w,h_{1}*})$ has a positive lower bound
uniformly. Combining~(\ref{s4.25}), (\ref{s4.26}), (\ref{s4.28}),
(\ref{s4.29}) and~(\ref{s4.34}) yields~(\ref{s4.31}).~~~$\square$\\

For convenient argument, we assume $s_2=s,~\gamma_1=\gamma_2=\gamma$ in the following Theorem.\\
\indent {\bf Theorem 4.4.}~~ Under the conditions of~Theorem 4.1, we further assume
that   $R(\Omega_{i})$ holds ($i=1,2,\cdots,l$), and
\begin{eqnarray}\label{s4.35}
w^{r}=\mathcal{O}(H^{r+s-1+\gamma}),~~~ h_{l}^{r+s-1}\gtrsim
H^{r+s-1+\gamma}.
\end{eqnarray}
Then there exists $u\in M(\lambda)$ and $u^{*}\in M^*(\lambda^*)$ such
that
\begin{eqnarray}\label{s4.36}
&&\|u^{w,h_{l}}-u\|_{1,\Omega}\lesssim h_{l}^{r+s-1},\\\label{s4.37}
&&\|u^{w,h_{l}}-u\|_{0,\Omega}\lesssim
H^{r+s-1+\gamma},\\\label{s4.38}
&&\|u^{w,h_{l}}-u\|_{1,\Omega\setminus\overline{F}}\lesssim
H^{r+s-1+\gamma},\\\label{s4.39}
&&\|u^{w,h_{l}*}-u^*\|_{1,\Omega}\lesssim
h_{l}^{r+s-1},\\\label{s4.40}
&&\|u^{w,h_{l}*}-u^*\|_{0,\Omega}\lesssim
H^{r+s-1+\gamma},\\\label{s4.41}
&&\|u^{w,h_{l}*}-u^*\|_{1,\Omega\setminus\overline{F}}\lesssim
H^{r+s-1+\gamma},\\\label{s4.42}
&&|\lambda^{w,h_{l}}-\lambda|\lesssim h_{l}^{2r+2s-2}.
\end{eqnarray}
\indent {\bf Proof.}~~Let~$u\in M(\lambda)$ and $u^*\in
M^*(\lambda^*)$,  such that~$u-u_{H}$ and $u^*-u_{H}^*$ both satisfy Lemma~2.2.
The proof of (\ref{s4.36})-(\ref{s4.42}) is completed by induction.\\
When $l=1$, Scheme 3.2 is actually Scheme 3.1. Hence, from
Theorem 4.1, Theorem 4.3 and (\ref{s4.35}) we know that (\ref{s4.36})-(\ref{s4.42}) hold for $l=0,1$.\\
Suppose (\ref{s4.36})-(\ref{s4.42}) hold for $l-2$, $l-1$, i.e.,
\begin{eqnarray}\label{s4.43}
&&\|u^{w,h_{l-2}}-u\|_{1,\Omega}\lesssim
h_{l-2}^{r+s-1},\\\label{s4.44}
&&\|u^{w,h_{l-2}}-u\|_{0,\Omega}\lesssim
H^{r+s-1+\gamma},\\\label{s4.45}
&&\|u^{w,h_{l-2}}-u\|_{1,\Omega\setminus\overline{F}}\lesssim
H^{r+s-1+\gamma},\\\label{s4.46}
&&\|u^{w,h_{l-2}*}-u^*\|_{1,\Omega}\lesssim
h_{l-2}^{r+s-1},\\\label{s4.47}
&&\|u^{w,h_{l-2}*}-u^*\|_{0,\Omega}\lesssim
H^{r+s-1+\gamma},\\\label{s4.48}
&&\|u^{w,h_{l-2}*}-u^*\|_{1,\Omega\setminus\overline{F}}\lesssim
H^{r+s-1+\gamma},\\\label{s4.49}
&&|\lambda^{w,h_{l-2}}-\lambda|\lesssim h_{l-2}^{2r+2s-2};
\end{eqnarray}
and
\begin{eqnarray}\label{s4.50}
&&\|u^{w,h_{l-1}}-u\|_{1,\Omega}\lesssim
h_{l-1}^{r+s-1},\\\label{s4.51}
&&\|u^{w,h_{l-1}}-u\|_{0,\Omega}\lesssim
H^{r+s-1+\gamma},\\\label{s4.52}
&&\|u^{w,h_{l-1}}-u\|_{1,\Omega\setminus\overline{F}}\lesssim
H^{r+s-1+\gamma},\\\label{s4.53}
&&\|u^{w,h_{l-1}*}-u^*\|_{1,\Omega}\lesssim
h_{l-1}^{r+s-1},\\\label{s4.54}
&&\|u^{w,h_{l-1}*}-u^*\|_{0,\Omega}\lesssim
H^{r+s-1+\gamma},\\\label{s4.55}
&&\|u^{w,h_{l-1}*}-u^*\|_{1,\Omega\setminus\overline{F}}\lesssim
H^{r+s-1+\gamma},\\\label{s4.56}
&&|\lambda^{w,h_{l-1}}-\lambda|\lesssim h_{l-1}^{2r+2s-2}.
\end{eqnarray}
Next we shall prove that (\ref{s4.36})-(\ref{s4.42}) hold for $l$.
Using the above formula and Lemma 2.4 to estimate the terms at the right hand side of (\ref{s4.9}) gives
\begin{eqnarray}\label{s4.57}
&&\|u^{w,h_{l}}-P_{h_{l}}u\|_{1,\Omega}\lesssim h_{l}^{r+s-1+\gamma}
+h_{l-1}^{\gamma}(h_{l}^{r+s-1}+h_{l-1}^{r+s-1})
+H^{r+s-1+\gamma}\nonumber\\
&&~~~~~~+H^{r+s-1+\gamma}+(H^{r+s-1+\gamma}+h_{l}^{r})
+H^{r+s-1+\gamma}\lesssim H^{r+s-1+\gamma}.
\end{eqnarray}
 Combining (\ref{s2.33}), (\ref{s2.34}) and (\ref{s4.57}) yields
~(\ref{s4.36}), (\ref{s4.37}) and (\ref{s4.38}).
By the same argument we can prove~(\ref{s4.39}), (\ref{s4.40}) and (\ref{s4.41}). From
(\ref{s2.26}), we have
\begin{eqnarray}\label{s4.58}
\lambda^{w,h_{l}}-\lambda=\frac{a(u^{w,h_{l}}-u,
u^{w,h_{l}*}-u^{*})}{b(u^{w,h_{l}},
u^{w,h_{l}*})}-\lambda\frac{b(u^{w,h_{l}}-u,
u^{w,h_{l}*}-u^{*})}{b(u^{w,h_{l}}, u^{w,h_{l}*})}.
\end{eqnarray}
Using the similar argument as that of Theorem 4.3 we know that
$b(u^{w,h_{l}}, u^{w,h_{l}*})$ has a positive lower bound
uniformly. Combing~(\ref{s4.36}), (\ref{s4.37}), (\ref{s4.39}),
(\ref{s4.40}) and~(\ref{s4.58}) yields~(\ref{s4.42}).~~~$\square$\\

{\bf Remark 4.1}.~~$\overline{\Omega}_{1}$ in (\ref{s3.7}) and
(\ref{s3.9}) can be different from that in (\ref{s3.6}) and
(\ref{s3.8}), which also ensure that
the corresponding estimates  in Theorem 4.2- Theorem 4.4 still hold.\\

{\bf Remark 4.2}.~~By dropping the steps computing $u_{H}^{*},
u^{w*}, e^{h*}, u^{w,h_{i}*},\lambda^{w,h_{i}*}$ in Scheme 3.2, and
by replacing $ \lambda^{w}=\frac{a(u^{w}, u^{w*})}{b(u^{w},
u^{w*})}$ and $\lambda^{w,h_{i}*}= \frac{a(u^{w,h_{i}},
u^{w,h_{i}*})}{a(u^{w,h_{i}}, u^{w,h_{i}*})}$ with $
\lambda^{w}=\frac{a(u^{w}, u^{w})}{b(u^{w}, u^{w})}$ and
$\lambda^{w,h_{i}}= \frac{a(u^{w,h_{i}},
u^{w,h_{i}})}{a(u^{w,h_{i}}, u^{w,h_{i}})}$, respectively, we are
able to establish multilevel discretizations based on local
defect-correction for symmetric Eigenvalue Problems.
Hence the corresponding estimates  in Theorem 4.1-Theorem 4.4 still hold.\\

{\bf Remark 4.3}.~~By referring to \cite{dai}, we can establish the
parallel version of Scheme 3.1 and Scheme 3.2 and have the
corresponding error estimates  in Theorem 4.3-Theorem 4.4
apparently.

\section{Numerical experiments}

\subsection{Computational method for $(\lambda_{H}^{*},u_{H}^{*})$ (see \cite{yang1,yang2})}

Assume that $(\lambda_{H},u_{H})$ is obtained from Scheme 3.1 or
Step 1 of Scheme 3.2, then $\lambda_{H}^{*}=\overline{\lambda}_{H}$,
and from \cite{yang1,yang2} we can  obtain $u_{H}^*$   by using the
following approach such that $|b(u_{H}, u_{H}^*)|$ has a positive
lower bound uniformly with respect to $H$.

Let $m_0$ be the algebraic multiplicity of
$\lambda_{H}$ and $l$ be ascent of $\lambda_{H}$.\\
\indent Let $u_{N}^{-}$ be the orthogonal projection of $u_{H}$ to
$N((\frac{1}{\lambda_{H}^{*}}-T_{H}^{*})^{l})$, and
$u_{H}^{*}=u_{H}^{-}/\|u_{H}^{-}\|_{0}$. When $u_{H}\in
N((\frac{1}{\lambda_{H}^{*}}-T_{H}^{*})^{l})$ it is clear that
$u_{H}^{*}=u_{H}$. When $u_H\notin
N((\frac{1}{\lambda_{H}^{*}}-T_{H}^{*})^{l})$, to find $u_{H}^{*}$,
first we seek a basis $\{\phi_i\}_{1}^{m_0}$ of
$N((\frac{1}{\lambda_{H}^{*}}-T_{H}^{*})^{l})$, and solve the
following equations
\begin{eqnarray}\label{s5.1}
\sum\limits_{i=1}^{m_0}\alpha_{i}(\phi_{i},\phi_{j})=(u_{H},\phi_{j}),~~~j=1,2,\cdots,m_0,
\end{eqnarray}
then let
\begin{eqnarray}\label{s5.2}
u_{H}^{-}&=&\sum\limits_{i=1}^{m_0}\alpha_{i}\phi_{i},\\\label{s5.3}
u_{H}^{*}&=&u_{H}^{-}/\|u_{H}^{-}\|_{0}.
\end{eqnarray}
Obviously, $u_{H}^{*}$ satisfies
\begin{eqnarray*}
(\frac{1}{\lambda_{H}^{*}}-T_{H}^{*})^{l}u_{H}^{*}=0,~~~\|u_{H}^{*}\|_{0}=1.
\end{eqnarray*}

Thus, to find $u_{H}^{*}$ which satisfying (\ref{s5.2}) and
(\ref{s5.3}) in $V_{H}^{0}(\Omega)$, the key is to seek a basis
$\{\phi_i\}_{1}^{m_0}$ of
$N((\frac{1}{\lambda_{H}^{*}}-T_{H}^{*})^{l})$. \\
\indent When $l=1$, it is
actually to solve the following equations to obtain a basis in the
solution space.
\begin{eqnarray}\label{s5.4}
&~&u_{H}^{(1)}\in V_{H}^{0}(\Omega),\nonumber\\
&~&a(v,u_{H}^{(1)})-\lambda_{H} b(v,u_{H}^{(1)})=0,~~~\forall v\in
V_{H}^{0}(\Omega).
\end{eqnarray}
(When $\lambda_{H}$ is a simple eigenvalue, $l=1$ and
$N((\frac{1}{\lambda_{H}^{*}}-T_{H}^{*})^l)$ is a one-dimensional
space spanned by the eigenfunction $u_{h}^{*}$.)\\
\indent When $l>1$, how to seek a basis $\{\phi\}_1^{m_0}$ of $N((\frac{1}{\lambda_{H}^{*}}-T_{H}^{*})^l)$ efficiently is an important issue of linear algebra.

\subsection{Numerical Examples}

\indent Consider the convection-diffusion equation
\begin{eqnarray}\label{s5.6}
-\Delta u+\mathbf{b}\cdot \nabla u=\lambda u,
~~~in~\Omega,~~~~u=0,~~~on~\partial\Omega,
\end{eqnarray}
where $\Omega=(-1,1)^2\setminus \{[0,1]\times[-1,0]\} $ or
$\Omega=(-1,1)^2\backslash\{\{0\}\times[-1,0]\}$. The first
eigenfunctions of both  problems have the singularities at the
origin. The exact eigenvalues, which are unknown, are thereby
replaced by approximate eigenvalues with high accuracy. For the
problem with $\mathbf{b}=(1, 1)^{T}$, $\mathbf{b}=(0, 3)^{T}$  and
$\mathbf{b}=(0, 10)^{T}$ on $\Omega=(-1,1)^2\setminus
\{[0,1]\times[-1,0]\} $, we take the approximate first eigenvalue as
$\lambda_{1}\approx 11.8897 $, $\lambda_{1}\approx 10.1397 $ and
$\lambda_{1}\approx34.6397 $, respectively. For the problem with
$\mathbf{b}=(1, 1)^{T}$, $\mathbf{b}=(0, 3)^{T}$ and $\mathbf{b}=(0,
10)^{T}$ on $\Omega=(-1,1)^2\backslash\{[0,1]\times\{0\}\}$, we take
the approximate first eigenvalue as $\lambda_{1}\approx10.621
$,$\lambda_{1}\approx 8.871 $ and $\lambda_{1}\approx 33.371 $,
respectively. We will report some numerical experiments by using
linear finite elements on uniform triangle meshes. In our numerical
experiments, we use Scheme 3.2 to solve the problem such that
$\Omega_{i}=(\frac{-1}{2^{i}},\frac{1}{2^{i}})^2\setminus
\{[0,\frac{1}{2^{i}}]\times[\frac{-1}{2^{i}},0]\} $ for L-shaped
domain,
$\Omega_{i}=(\frac{-1}{2^{i}},\frac{1}{2^{i}})^2\backslash\{\{0\}\times[-\frac{1}{2^{i}},0]\}$
for slit domain, $i=1,2,\cdots,6$,
and locally fine grids have the same degree of freedom as that of globally mesoscopic grid (see Tables 1-6).\\
\indent In our experiments, according to the assumptions of Theorem
4.4, we approximately take $\gamma_1=\gamma_2=1/2,2/3$ and
$s=s_2=1/2,2/3$ so that (\ref{s4.35}) holds for slit domain and
L-shaped domain, respectively.   We use MATLAB 2011b under the
package of Chen (see \cite{chen}) to solve the problems, and the
numerical results are shown in Tables 1-6. From Tables 1-4 we can see
that without increasing degree of freedom on locally fine grids, the
first local defect correction can largely improve the accuracy of
the eigenvalues, and the local defect corrections that follows  can
gradually improve the accuracy of the eigenvalues by overcoming the
singularity at the origin. But Tables 5-6 also indicate that Scheme
3.2 is not valid for the problems with $\mathbf{b}=(0,10)^T$.
Concerning this point,  the figures of eigenfunction and its adjoint
pair (see Fig. 5.1-5.2) shows their  functional-value abrupt changes
mainly center on boundary layer, which may lead to the invalidity of
Scheme 3.2.
 for the case  $\mathbf{b}=(0,10)^T$, we adopt the parallel version of Scheme 3.2 to make local defect-corrections on boundary layers with functional-value abrupt changes (see also Fig. 5.1-5.2). \\\indent
 Specifically speaking,  for the L-shaped domain, we find that it's better to make local defect-corrections near the origin on slightly small area $\Omega^1_{i}=(\frac{-1}{2^{i+1}},\frac{1}{2^{i+1}})^2\backslash\{\{0\}\times[-\frac{1}{2^{i+1}},0]\}$
 for both eigenfunction and its adjoint eigenfunction;
 as for the other  local defect-correction areas, we set as
$\Omega^2_{i}=(\frac{-1}{2^{i-1}},\frac{1}{2^{i-1}})\times(1-\frac{1}{2^{i}},1)$ for the eigenfunction,
 $\Omega^3_{i}=(-\frac{1}{2}-\frac{1}{2^{i} },-\frac{1}{2}+\frac{1}{2^{i}})\times(-1,-1+\frac{1}{2^{i}})$ for the adjoint eigenfunction,
 respectively; the related numerical results is given in Table 7.
 Here we set
 \begin{eqnarray*}
  DOF_w \approx \frac{3}{4}\times DOF_{\Omega^{1}_i}\approx \frac{3}{2} \times DOF_{\Omega^{2}_i}\approx 4\times DOF_{\Omega^{3}_i}~(i=1,2,\cdots).
 \end{eqnarray*}
 \indent
 For the slit domain, we set as the local defect-correction area $\Omega^1_{i}=(\frac{-1}{2^{i}},\frac{1}{2^{i}})^2\backslash\{\{0\}\times[-\frac{1}{2^{i}},0]\}$
 for both eigenfunction and its adjoint eigenfunction,
$\Omega^2_{i}=(\frac{-1}{2^{i-1}},\frac{1}{2^{i-1}})\times(1-\frac{1}{2^{i}},1)$ for the eigenfunction,
 $\Omega^3_{i}=(\frac{1}{2}-\frac{1}{2^{i} },\frac{1}{2}+\frac{1}{2^{i}})\times(-1,-1+\frac{1}{2^{i}})$ and
 $\Omega^4_{i}=(-\frac{1}{2}-\frac{1}{2^{i} },-\frac{1}{2}+\frac{1}{2^{i}})\times(-1,-1+\frac{1}{2^{i}})$ for the adjoint eigenfunction,
 respectively; the related numerical results is given in Table 8.  Here we set
 \begin{eqnarray*}
  DOF_w=DOF_{\Omega^{1}_i}=DOF_{\Omega^{2}_i}\approx 2\times DOF_{\Omega^{3}_i}=2\times DOF_{\Omega^{4}_i}~(i=1,2,\cdots).
 \end{eqnarray*}
 \indent Table 7.12 and Table 7.16 in \cite{carstensen}  show that, using the adaptive homotopy method to solve the L-shaped domain problem with $\mathbf{b}=(10,0)^T$,
 the approximate eigenvalue can have 4-5 significant digits with $DOF=154994$ and $124469$, thus the adaptive homotopy method is efficient. However, by using our algorithm,
 the approximate eigenvalue can have 6 significant digits with $DOF_H=12033$ (see Table 7), which also indicates  our algorithm is  efficient.
\\
\\
\begin{tabular}{lllllllll}
\multicolumn{6}{c}{\bf Table 1: {   $\Omega=(-1,1)^2\backslash\{[0,1]\times[-1,0]\}$, $b=(0,3)^T$.}}\\
\hline\noalign{\smallskip}
$DOF_H$&$DOF_w$&$\lambda_H$&$\lambda^w$&$\lambda^{w,h_1}$\\
\noalign{\smallskip}\hline\noalign{\smallskip}
705&2945&11.94916&11.91247&11.89949\\
2945&12033&11.91250&11.89859&11.89343\\
12033&195585& 11.89859&11.89109&11.89028\\
\noalign{\smallskip}\hline\noalign{\smallskip}
$\lambda^{w,h_2}$&$\lambda^{w,h_3}$&$\lambda^{w,h_4}$&$\lambda^{w,h_5}$&$\lambda^{w,h_6}$\\
\noalign{\smallskip}\hline\noalign{\smallskip}
 11.89466&11.89275&-&-&-\\
 11.89146& 11.89068&11.89037&-&-\\
 11.88996& 11.88983& 11.88978&-&-\\
\noalign{\smallskip}\hline
\end{tabular}
\\\\\\
\begin{tabular}{llllllllll}
\multicolumn{6}{c}{\bf Table 2: {   $\Omega=(-1,1)^2\backslash\{[0,1]\times[-1,0]\}$, $b=(1,1)^T$.}}\\
\hline\noalign{\smallskip}
$DOF_H$&$DOF_w$&$\lambda_H$&$\lambda^w$&$\lambda^{w,h_1}$\\
\noalign{\smallskip}\hline\noalign{\smallskip}
705&2945&10.21836&10.16730&10.15332\\
2945&12033&10.16730&10.14979&10.14438\\
12033&195585& 10.14979&10.14117&10.14034\\
\hline\noalign{\smallskip}
$\lambda^{w,h_2}$&$\lambda^{w,h_3}$&$\lambda^{w,h_4}$&$\lambda^{w,h_5}$&$\lambda^{w,h_6}$\\
\noalign{\smallskip}\hline\noalign{\smallskip}
 10.14836&10.14643&-&-&-\\
 10.14238&10.14160&10.14129&-&-\\
 10.14002&10.13989&10.13984&-&-\\
\noalign{\smallskip}\hline
\end{tabular}
\\\\\\
\begin{tabular}{lllllllll}%
\multicolumn{5}{c}{\bf Table 3: {   $\Omega=(-1,1)^2\backslash\{\{0\}\times[-1,0]\}$, $b=(0,3)^T$.}}\\
\hline\noalign{\smallskip}
$DOF_H$&$DOF_w$&$\lambda_H$&$\lambda^w$&$\lambda^{w,h_1}$\\
\noalign{\smallskip}\hline\noalign{\smallskip}
 945&3937&10.79397&10.70640&10.66393\\
 3937&16065&10.70630&10.66356&10.64247\\
 16065&64897&10.66353&10.64237&10.63186\\
\noalign{\smallskip}\hline\noalign{\smallskip}
$\lambda^{w,h_2}$&$\lambda^{w,h_3}$&$\lambda^{w,h_4}$&$\lambda^{w,h_5}$&$\lambda^{w,h_6}$\\
\noalign{\smallskip}\hline\noalign{\smallskip}
10.64333&  10.63315&-&-&-\\
10.63208& 10.62691&10.62434&-&-\\
10.62664&10.62404&10.62274&10.62209&-\\
\noalign{\smallskip}\hline
\end{tabular}
\\\\\\
\begin{tabular}{llllllll}%
\multicolumn{6}{c}{\bf Table 4: {   $\Omega=(-1,1)^2\backslash\{\{0\}\times[-1,0]\}$, $b=(1,1)^T$.}}\\
\hline\noalign{\smallskip}
$DOF_H$&$DOF_w$&$\lambda_H$&$\lambda^w$&$\lambda^{w,h_1}$\\
\noalign{\smallskip}\hline\noalign{\smallskip}
 945&3937&9.06244&8.96099& 8.91741\\
 3937&16065& 8.96090&8.91470&8.89333\\
 16065&64897& 8.91468&8.89266&8.88207\\
\noalign{\smallskip}\hline\noalign{\smallskip}
$\lambda^{w,h_2}$&$\lambda^{w,h_3}$&$\lambda^{w,h_4}$&$\lambda^{w,h_5}$&$\lambda^{w,h_6}$\\
\noalign{\smallskip}\hline\noalign{\smallskip}
8.89663& 8.88641&-&-&-\\
8.88289&8.87772&8.87514&-&-\\
8.87684& 8.87424&8.87294& 8.87229&-\\
\noalign{\smallskip}\hline
\end{tabular}\\\\\\
\begin{tabular}{llllllllll}%
\multicolumn{6}{c}{\bf Table 5: {   $\Omega=(-1,1)^2\backslash\{[0,1]\times[-1,0]\}$, $b=(0,10)^T$. }}\\
\hline\noalign{\smallskip}
$DOF_H$&$DOF_w$&$\lambda_H$&$\lambda^w$&$\lambda^{w,h_1}$\\
\noalign{\smallskip}\hline\noalign{\smallskip}
705&2945&34.58756&34.63484&34.61999\\
2945&12033&34.63473&34.64168&34.63605\\
12033&195585& 34.64167&34.64066 &34.63982\\
\noalign{\smallskip}\hline\noalign{\smallskip}
$\lambda^{w,h_2}$&$\lambda^{w,h_3}$&$\lambda^{w,h_4}$&$\lambda^{w,h_5}$&$\lambda^{w,h_6}$\\
\noalign{\smallskip}\hline\noalign{\smallskip}
 34.61632& 34.61460&-&-&-\\
 34.63437&34.63363&34.63333&-&-\\
 34.63952&34.63939&34.63934&-&-\\
\noalign{\smallskip}\hline
\end{tabular}
\\
\\\\\\
\begin{tabular}{lllllllll}%
\multicolumn{6}{c}{\bf Table 6: {   $\Omega=(-1,1)^2\backslash\{\{0\}\times[-1,0]\}$, $b=(0,10)^T$. }}\\
\hline\noalign{\smallskip}
$DOF_H$&$DOF_w$&$\lambda_H$&$\lambda^w$&$\lambda^{w,h_1}$\\
\noalign{\smallskip}\hline\noalign{\smallskip}
 945&3937&  33.43287&33.42950&33.38763\\
 3937&16065&33.42885&33.40686&33.38587\\
 16065&64897&33.40671& 33.39070&33.38021\\
\hline\noalign{\smallskip}
$\lambda^{w,h_2}$&$\lambda^{w,h_3}$&$\lambda^{w,h_4}$&$\lambda^{w,h_5}$&$\lambda^{w,h_6}$\\
\noalign{\smallskip}\hline\noalign{\smallskip}
33.36894&  33.35923&-&-&-\\
33.37595& 33.37090&33.36836&-&-\\
33.37510& 33.37253&33.37124& 33.3706&-\\
\noalign{\smallskip}\hline
\end{tabular}\\\\\\
\begin{tabular}{lllllllllll}%
\multicolumn{6}{c}{\bf Table 7: {   $\Omega=(-1,1)^2\backslash\{[0,1]\times[-1,0]\}$, $b=(0,10)^T$.}}\\
\hline\noalign{\smallskip}
$DOF_H$&$DOF_w$& $\lambda_H$&$\lambda^w$&$\lambda^{w,h_1}$\\
\noalign{\smallskip}\hline\noalign{\smallskip}
705&2945&34.58756&34.63484&34.64245\\
2945&12033& 34.63474&34.64169&34.64162\\
12033&195585& 34.64166&34.64067&34.64017\\
\noalign{\smallskip}\hline\noalign{\smallskip}
$\lambda^{w,h_2}$&$\lambda^{w,h_3}$&$\lambda^{w,h_4}$&$\lambda^{w,h_5}$&$\lambda^{w,h_6}$\\
\noalign{\smallskip}\hline\noalign{\smallskip}
  34.64102& 34.63949&-&-&-\\
 34.64049& 34.63980& 34.63951&-&-\\
  34.63990&34.63978&34.63973&-&- \\
\noalign{\smallskip}\hline
\end{tabular}\\\\\\
\begin{tabular}{lllllllllll}%
\multicolumn{5}{c}{\bf Table 8: {   $\Omega=(-1,1)^2\backslash\{\{0\}\times[-1,0]\}$, $b=(0,10)^T$.}}\\
\hline\noalign{\smallskip}
$DOF_H$&$DOF_w$ &$\lambda_H$&$\lambda^w$&$\lambda^{w,h_1}$\\
\noalign{\smallskip}\hline\noalign{\smallskip}
 945&3937& 33.43287&33.42950&33.40090\\
 3937&16065& 33.42885&33.40686&33.38916\\
16065&64897&33.40671& 33.39070&33.38103\\
\noalign{\smallskip}\hline\noalign{\smallskip}
$\lambda^{w,h_2}$&$\lambda^{w,h_3}$&$\lambda^{w,h_4}$&$\lambda^{w,h_5}$&$\lambda^{w,h_6}$\\
\noalign{\smallskip}\hline\noalign{\smallskip}
 33.38429&33.37474&-&-&-\\
33.37976&33.37475&33.37221&-&-\\
33.37605&33.37349& 33.37220& 33.37155&-\\
\noalign{\smallskip}\hline
\end{tabular}\\\\\\
\begin{tabular}{cc}
\includegraphics[width=0.4\textwidth]{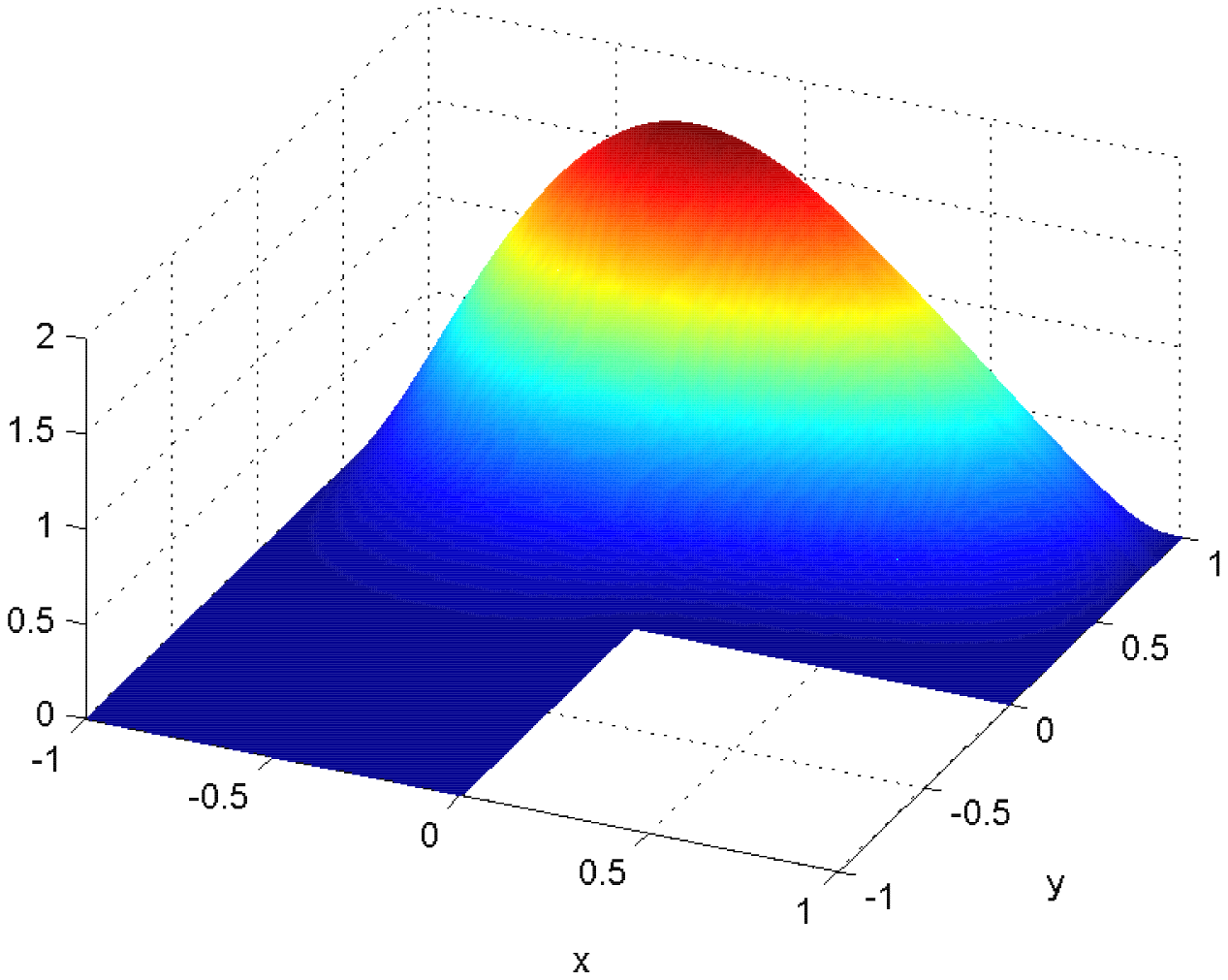}&
\includegraphics[width=0.4\textwidth]{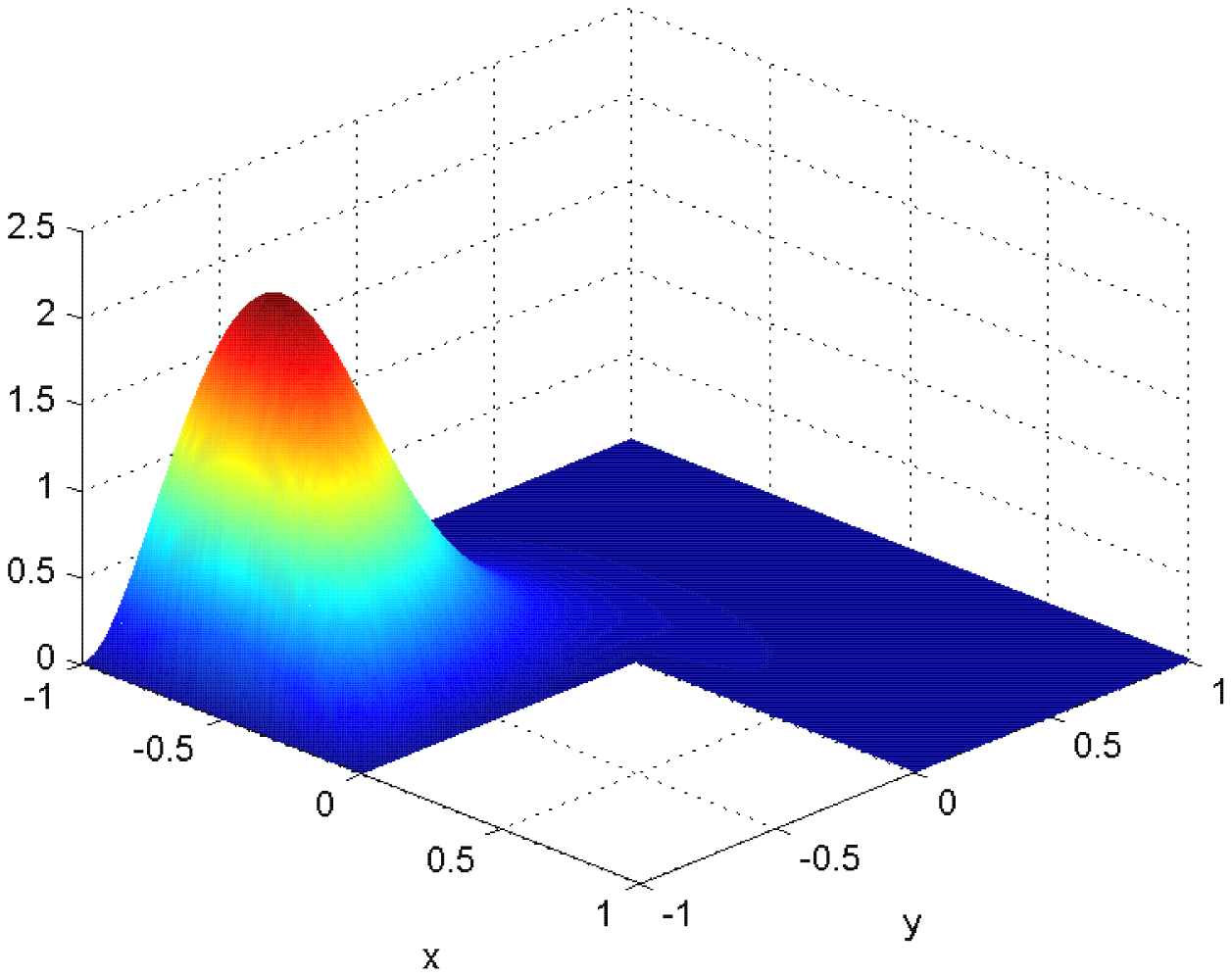}
\end{tabular}
\\{\textrm{\small\bf Fig. 5.1.  Eigenfuntion and its adjoint pair with $\mathbf{b}=(0,10)^T$ on $\Omega=(-1,1)^2\backslash
\{[0,1]\times[-1,0]\}$}}
\\
\begin{tabular}{cc}
\includegraphics[width=0.4\textwidth]{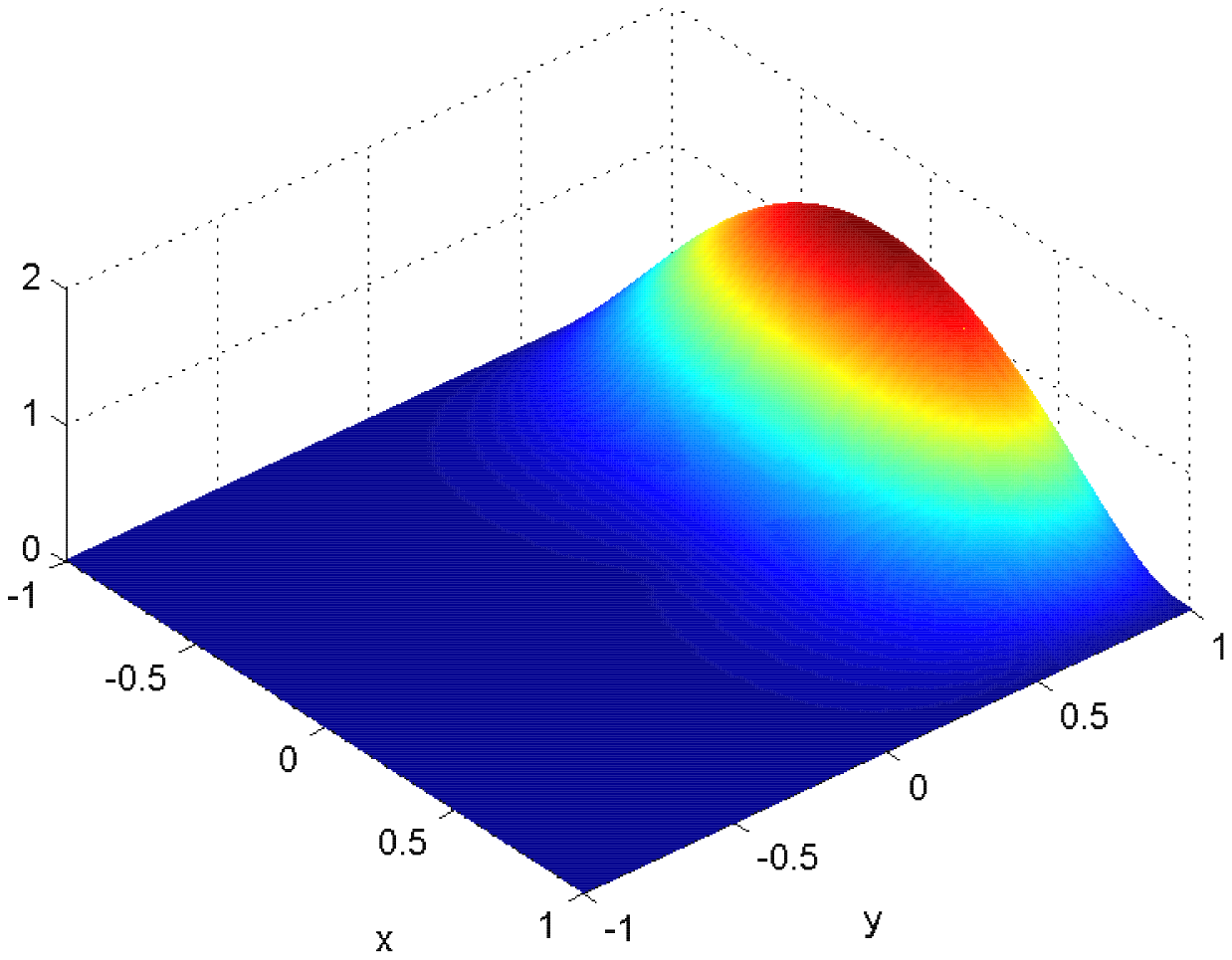}&
\includegraphics[width=0.4\textwidth]{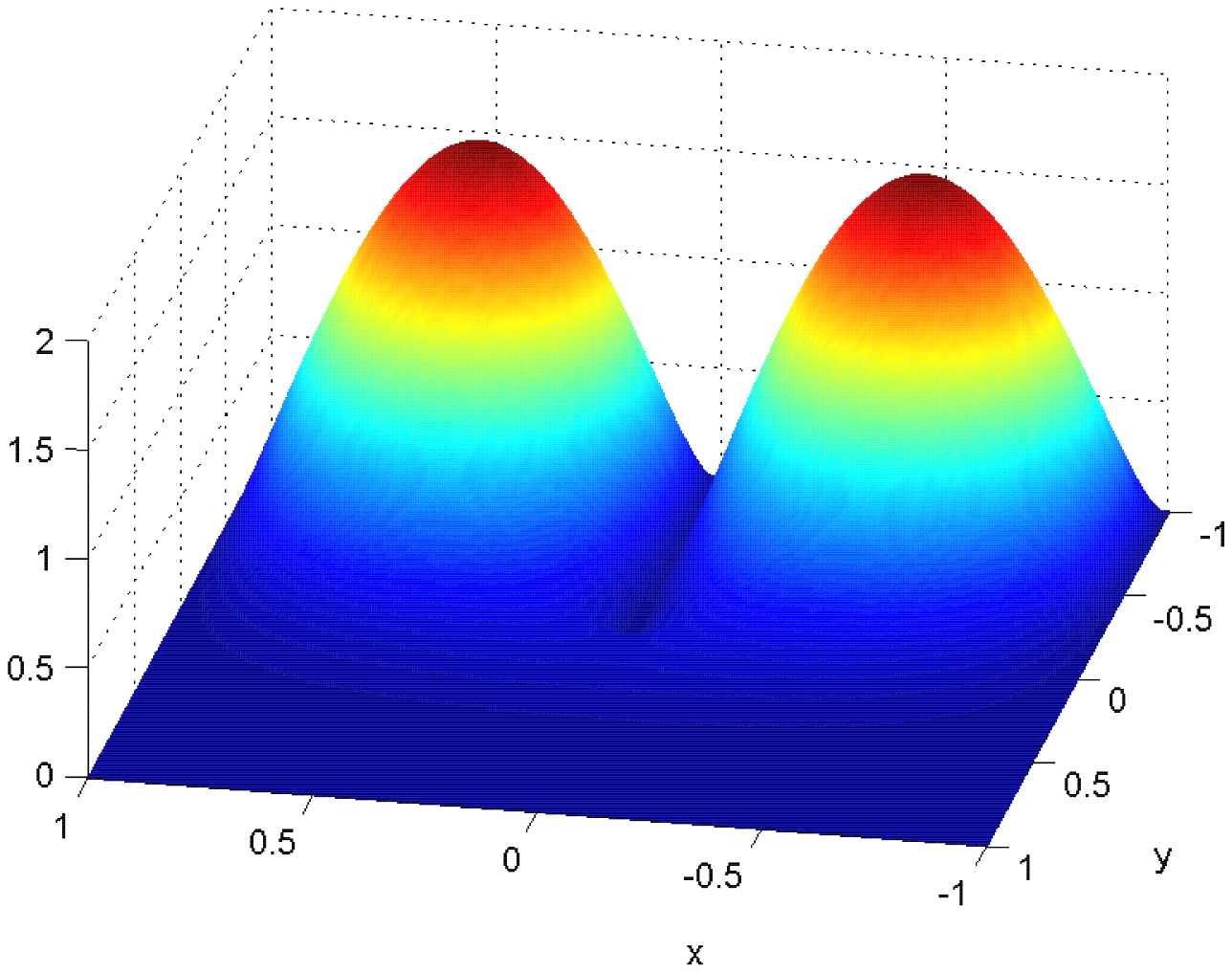}
\end{tabular}
\\{\textrm{\small\bf Fig. 5.2.  Eigenfuntion and its adjoint pair with $\mathbf{b}=(0,10)^T$ on $\Omega=(-1,1)^2\backslash\{\{0\}\times[-1,0]\}$}}

\end{document}